\documentclass[british,11pt]{article}

\usepackage[english,british]{babel}
\usepackage{amsmath,amssymb,enumerate,amsthm}
\usepackage{cite}
\usepackage[latin1]{inputenc}
\usepackage[T1]{fontenc}

\usepackage{psfig,labelfig}
\usepackage{graphicx}
\usepackage[dvips]{epsfig}

\newtheorem{thm}{Theorem}[section]

\newtheorem{lem}[thm]{Lemma}
\newtheorem{defi}[thm]{Definition}
\newtheorem{exa}[thm]{Example}
\newtheorem*{cla}{Claim}

\newcommand{\C}{{\mathbb C}}
\newcommand{\Hy}{{\mathbb H}}
\newcommand{\Sp}{{\mathbb S}}
\newcommand{\R}{{\mathbb R}}

\newcommand{\Z}{{\mathbb Z}}
\newcommand{\N}{{\mathbb N}}

\DeclareMathOperator{\arcsinh}{arcsinh}

\DeclareMathOperator{\arctanh}{arctanh}

\DeclareMathOperator{\ir}{Int}
\DeclareMathOperator{\lip}{Lip}
\DeclareMathOperator{\ao}{\alpha_1}
\DeclareMathOperator{\at}{\alpha_2}
\DeclareMathOperator{\dist}{dist}

\DeclareMathOperator{\capa}{cap}

\DeclareMathOperator{\iso}{Isom}

\title{The Jacobian of a Riemann surface and the geometry of the cut locus of simple closed geodesics}

\author{Bjoern Muetzel \thanks{E-mail address : bjorn.mutzel@gmail.com}\\
 \\
\small Department of Mathematics,  Dartmouth College, Hanover, NH 03755, USA\\[-0.8ex]}

\begin{document}
\maketitle

\begin{abstract}
To any compact Riemann surface of genus $g$ one may assign a principally polarized abelian variety of dimension $g$, the \textit{Jacobian} of the Riemann surface. The Jacobian is a complex torus, and a Gram matrix of the lattice of a Jacobian is called a \textit{period Gram matrix}. This paper provides upper and lower bounds for all the entries of the period Gram matrix with respect to a suitable homology basis. These bounds depend on the geometry of the cut locus of non-separating simple closed geodesics. Assuming that the cut loci can be calculated, a theoretical approach is presented followed by an example where the upper bound is sharp. Finally we give practical estimates based on the Fenchel-Nielsen coordinates of surfaces of signature $(1,1)$, or \textit{Q-pieces}. The methods developed here have been applied to surfaces that contain small non-separating simple closed geodesics in \cite{bmm}.
\end{abstract}

\section{Introduction}

Let $S$ be a hyperbolic Riemann surface of genus $g \geq 2$. We call a set of $2g$ oriented simple closed geodesics
\[
 {\rm A} = (\alpha_1, \alpha_{2},...,\alpha_{2g-1},\alpha_{2g})
\]
a \textit{canonical basis}, if
\begin{itemize}
\item[-] for each $\alpha_i$ there exists exactly one $\alpha _{\tau(i)}=
\left\{ {\begin{array}{*{20}c}
   {\alpha_{i+1}}  \\
   {\alpha_{i-1}}  \\
\end{array}} \right. \text{  if  } \begin{array}{*{20}c}
   {i \text{  odd } }  \\
   {i \text{  even }}  \\
\end{array}
 \in  {\rm A}$ that intersects $\alpha_i$ in exactly one point.
\item[-] the curves are oriented in a way, such that
\[
\ir ( \alpha_{i}, \alpha_{i+1}) = 1  \text{ \ for all \ } i= 1,3,...,2g-1,
\]
where $\ir(\cdot,\cdot)$ denotes the algebraic intersection number.
\end{itemize}

Note that {\rm A} can be called a basis as the homology classes $\left( [\alpha _i] \right)_{i = 1,...,2g} \subset H_1(S,\Z)$ form a basis of $H_1(S,\R)$ as a vector space. In the vector space of real harmonic 1-forms on $S$, let $\left( {\sigma _k } \right)_{k = 1,...,2g}$ be the \textit{dual basis for} $\left( {[\alpha _i] } \right)_{i = 1,...,2g} \subset H_1(S,\Z)$ defined by
\[
\int\limits_{[\alpha_i]} {\sigma_k } = \delta_{ik} .
\]
A \textit{period Gram matrix} $P_S$ (with respect to ${\rm A}$) of $S$ is the Gram matrix
\begin{equation}
P_S=\left(\left\langle {\sigma_i,\sigma_j} \right\rangle  \right)_{i,j= 1,...,2g}=\left( {\int\limits_S {\sigma_i  \wedge {}^ * } \sigma_j } \right)_{i,j= 1,...,2g}.   \nonumber
\end{equation}

\begin{figure}[h!]
\SetLabels
\L(.30*.90) $\mathcal{Q}_i$\\
\L(.66*.54) $\alpha_i$\\
\L(.43*.88) $\,\,\alpha_{\tau(i)}$\\
\L(.61*.80) $\,\alpha_{i\tau(i)}$\\
\L(.60*.03) $\beta_i$\\
\endSetLabels
\AffixLabels{%
\centerline{%
\includegraphics[height=6cm,width=8cm]{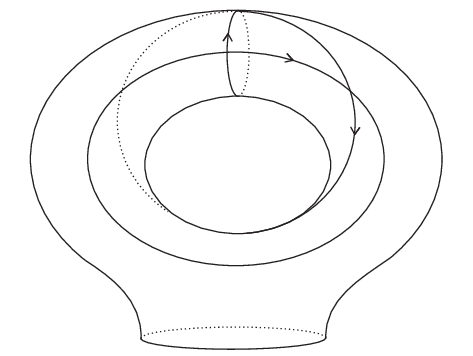}}}
\caption{A Q-piece $\mathcal{Q}_i$. The curve $\alpha_{i\tau(i)}$ is in the free homotopy class of $\alpha_i(\alpha_{\tau(i)})^{-1}$.}
\label{fig:Qj}
\end{figure}
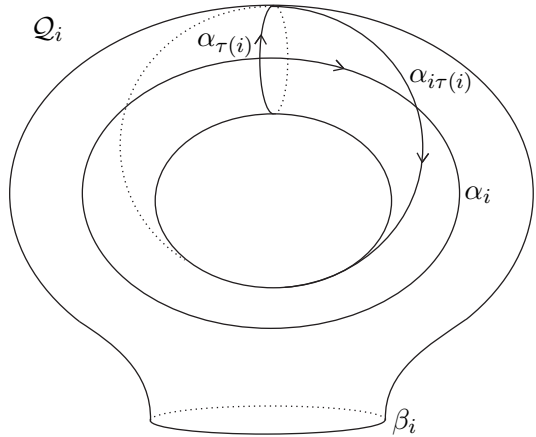

This period matrix $P_S$ defines a complex torus, the \textit{Jacobian} or \textit{Jacobian variety} $J(S)$ of the Riemann surface $S$ (see \cite{fk}, chapter III). Let
\[
E(\sigma_i)= E_S(\sigma_i)=\int\limits_S {\sigma_i  \wedge {}^ * } \sigma_i = \left\langle {\sigma_i,\sigma_i} \right\rangle
\]
be the \textit{energy of $\sigma_i$ (over $S$)}. As $P_S$ is a Gram matrix, $E(\sigma_i)$ is also the squared norm of a vector $v_i$ in the lattice of the Jacobian.\\
In this paper, we examine the connection between the metric, hyperbolic geometry of a compact Riemann surface and the geometry of its Jacobian. In previous papers (see \cite{bsi} or \cite{se}), this approach has been taken for special cases, for example when the Riemann surface is a real algebraic curve. For these special cases, there exist algorithms to calculate the period matrix.\\
 On the contrary, when the surface is given in terms of its hyperbolic metric we do not know explicitly the harmonic 1-forms. As the exact computation of the period Gram matrix seems very difficult, we look for an approximation or estimate. In this paper we find upper and lower bounds for all entries of the period Gram matrix based on the hyperbolic metric of an arbitrary compact Riemann surface. The bounds depend on the geometry of the cut loci of the curves in a canonical basis and related simple closed geodesics. They are obtained by estimating the energy of the corresponding dual harmonic forms.\\ 
The paper is divided into four parts: Section 1 is the introduction and Section 2 contains the preliminaries. In Section 3 we give theoretical estimates on the entries of the period Gram matrix. Here theoretical means that the estimates depend on some geometrical quantities that are not explicitly given in Fenchel-Nielsen coordinates. Finally, in Section 4 we apply the approach of the previous section to find explicit estimates.\\
More precisely, in Section 3 we find upper bounds for the energy of the dual harmonic forms by estimating the capacity of hyperbolic tubes as follows: let $T(\alpha_{\tau(i)}) \subset S$ be a topological tube, embedded in $S$ that contains the geodesic $\alpha_{\tau(i)}$ in its interior. The capacity of such a tube gives an upper bound for the energy $E(\sigma_i)$ of $\sigma_i$. This is the diagonal entry $p_{ii}$ of the period Gram matrix $P_S$:
\[
       \capa(T(\alpha_{\tau(i)})) \geq E(\sigma_i)=p_{ii}.
\]
In our theoretical approach, the boundary of such a tube will be provided by the cut locus of a simple closed geodesic of the canonical basis. More precisely, we will take $T(\alpha_{\tau(i)})=S_{\tau(i)}$, where $S_{\tau(i)}$ is the surface obtained by cutting open $S$ along the cut locus $CL(\alpha_{\tau(i)})$ of $\alpha_{\tau(i)}$ (see (\ref{eq:cut_locus})). This allows us to extend our tubes over the whole surface $S$ and to obtain a lower bound on $E(\sigma_i)$. This bound is obtained using projections of vector fields onto curves. Upper and lower bounds for the non-diagonal elements are obtained in a similar way with the help of the polarization identity.\\
The method presented in Section 3 relies on the premise that the cut loci in question can be calculated. The quality of the approximation depends on the geometry of the surface. This is illustrated by two examples. One based on a \textit{necklace surface} and one based on a \textit{linear surface}  presented in this section. We obtain the following estimates which - to our knowledge - are new in the literature. Here \textbf{Example 3.1} shows the limitations of the method, while \textbf{Example 3.2} shows a case where the upper bound is sharp.
\begin{thm}
Let $N$ be a necklace surface of genus $g \geq 2$. Then there is a canonical basis  ${\rm A} = (\alpha_i)_{i=1,...,2g}$ for which we have: if $N_1$ is the surface obtained by cutting open $N$ along the cut locus $CL(\alpha_1)$ of $\alpha_1$ and $P_{N}= (p_{ij})_{i,j}$ is the period Gram matrix with respect to ${\rm A}$. Then
\[
    \frac{c_{\alpha_1}}{g-1}   \geq  p_{22} \geq \frac{\|[\alpha_1]\|^2_s}{4\pi(g-1)} \text{ \ \ and \ \ } \capa(N_1) \geq p_{22}, \text{ \ \ but \ \ }   \capa(N_1) \geq \frac{\ell(\alpha_1)}{\pi},
\]
where $c_{\alpha_1}$ is a factor that depends only on the fixed length $\ell(\ao)$ of $\alpha_1$ and $\|[\alpha_1]\|_s$ is the length of a shortest multicurve in the same homology class as $\alpha_1$.
\label{thm:necklace}
\end{thm}
Hence $p_{22}$ is of order $\frac{1}{g}$ and goes to zero, as $g$ goes to infinity. Our upper bound, on the contrary, is always bigger than the constant $\frac{\ell(\alpha_1)}{\pi}$. This example shows an instance of the case where our upper bound cannot be of the right order.
\begin{thm}
Let $L$ be a linear surface of genus $g \geq 2$. Then there is a canonical basis  ${\rm A} = (\alpha_i)_{i=1,...,2g}$ for which we have: if $L_1$ is the surface obtained by cutting open $L$ along the cut locus $CL(\alpha_1)$ of $\alpha_1$ and $P_{L}= (p_{ij})_{i,j}$ is the period Gram matrix with respect to ${\rm A}$. Then 
\[
     p_{22} = \capa(L_1) - \epsilon_L,
\]
where $\epsilon_L > 0$ depends on the geometry of $L$ and may become arbitrarily small.
\label{thm:linear}
\end{thm}
This example shows an instance of the case where the bound is sharp for any genus.\\
The methods developed in this paper have been applied to surfaces that contain a short simple closed geodesic $\gamma$ in \cite{bmm}. If $\gamma$ is a separating closed geodesic then the matrix $P_S$ converges to a block matrix if $\ell(\gamma)$ goes to zero. In this case the bound on a non-diagonal entry of $P_S$ given with respect to a suitable canonical basis is sharp. The results of this section are summarized in \textbf{Theorem \ref{thm:diag}} and \textbf{\ref{thm:ndiag}}.\\
Finally, in Section 4 we apply our approach to a surface $S$ given in Fenchel-Nielsen coordinates. The theoretical estimates from Section 3 depend on the cut loci of curves which are often difficult to handle. To bypass this problem we work with one-holed tori or \textit{Q-pieces}. Under this condition the cut loci of the elements of a canonical basis can be (at least partially) calculated.\\
The method does not use all $6g-6$ Fenchel-Nielsen coordinates, but only uses the $3g$ coordinates of $g$ Q-pieces which are determined by two intersecting homology classes. 
More precisely, let
\[
(\mathcal{Q}_i)_{i=1,3,...,2g-1}\subset S
\]
be a set of Q-pieces, whose interiors are pairwise disjoint.
Let $\beta_i$ be the boundary geodesic of $\mathcal{Q}_i$, $\alpha_i$ an interior simple closed geodesic, and ${\bf tw}_i$ the twist parameter at $\alpha_i$. The geometry of $\mathcal{Q}_i$ is determined by the Fenchel-Nielsen coordinates $(\ell(\beta_i),\ell(\alpha_i),{\bf tw}_i)$. We assume furthermore that
\begin{equation}
\cosh(\frac{\ell(\alpha_i)}{2}) \leq \cosh(\frac{\ell(\beta_i)}{6}) + \frac{1}{2} \text{ \ for all \ } i \in  \{1,3,...,2g-1\}.
\label{eq:parlier1}
\end{equation}
By \cite{sch}, \textbf{Corollary 4.1} such a pair $(\alpha_i,\beta_i)$ always exists. This choice is made for two reasons. First it facilitates the calculation of the length of a suitable $\alpha_{\tau(i)}=\alpha_{i+1}$. Second it follows from the collar lemma in hyperbolic geometry that small simple closed geodesics have large collars, which in return gives good estimates for our upper bounds on the energies (see Section 2.2). \\
In Section 4, we first determine suitable $\alpha_{\tau(i)} \subset \mathcal{Q}_i$ for each $\alpha_i$, such that the pairs\\ $((\alpha_i,\alpha_{\tau(i)}))_{i=1,3,...,2g-1}$ form a canonical basis. Now fix an $i \in  \{1,3,...,2g-1\}.$ Let $\alpha_{i\tau(i)} \subset \mathcal{Q}_i$ be the simple closed geodesic in the free homotopy class of $\alpha_i(\alpha_{\tau(i)})^{-1}$ (see \textit{Figure \ref{fig:Qj}}). For $j \in \{i,\tau(i),i\tau(i)\}$, let
\begin{itemize}
\item[-]  $\beta_j=\beta_i$ be the boundary geodesic of $\mathcal{Q}_i$,
\item[-]  ${\bf tw}_j$ the twist parameter at $\alpha_j$,
\item[-]  $FN_j := (\ell(\beta_j),\ell(\alpha_j),{\bf tw}_j)$ the corresponding Fenchel-Nielsen coordinates of $\mathcal{Q}_i$.
\end{itemize}

In Section 4.1,  $FN_{\tau(i)}$ and $FN_{i\tau(i)}$ are calculated from $FN_{i}$. Section 4.2 and 4.3 give explicit functions providing upper and lower bounds on all entries of $P_S=(p_{ij})_{i,j}$. These estimates are summarized in \textbf{Theorem \ref{thm:period_Q}}.\\
An example of a period Gram matrix obtained via this method is worked out in \textbf{Example 4.3}. We obtain the following estimate for a parameter family of Riemann surfaces of genus $2$.
\begin{exa}
Let $\mathcal{Q}_1$ and $\mathcal{Q}_3$ be two isometric Q-pieces given in Fenchel-Nielsen coordinates $FN_1$ and $FN_3$, respectively, where $FN_i = (\ell(\beta_i),\ell(\alpha_i),{\bf tw}_i) = (2,1,0.1)$ for  $i\in \{1,3\}$.
Let $S=\mathcal{Q}_1+\mathcal{Q}_3$ be a Riemann surface of genus $2$, which we obtain by gluing $\mathcal{Q}_1$ and $\mathcal{Q}_3$ along $\beta_1$ and $\beta_3$ with arbitrary twist parameter ${\bf tw}_{\beta} \in (-\frac{1}{2},\frac{1}{2}]$. Then there exists a canonical basis ${\rm A}= (\alpha_1,...,\alpha_4)$ and a corresponding period Gram matrix $P_S$, such that
\[
\left( {\begin{array}{*{20}c}
 2.11 & -0.46 & -0.42 & -0.26 \\
-0.46 & 0.33 & -0.26 & -0.11 \\
-0.42 & -0.26 & 2.11 & -0.46 \\
-0.26 & -0.11 & -0.46 & 0.33 \\
\end{array}} \right)
\leq P_S \leq
\left( {\begin{array}{*{20}c}
2.53 & 0.20 & 0.42 & 0.26 \\
0.20 & 0.44 & 0.26 & 0.11 \\
0.42 & 0.26 & 2.53 & 0.20 \\
0.26 & 0.11 & 0.20 & 0.44 \\
\end{array}} \right).
\]
\end{exa}
In Section 4.4, the results are summarized in \textit{Table \ref{tab:notwist}} and compared with the upper bound that can be obtained from the method in \cite{bs} applied to Q-pieces. This bound is in general much larger.\\
We note that for $j \in \{i,i+1\}$ our upper bound for the diagonal entry $p_{\tau(j)\tau(j)} =E(\sigma_{\tau(j)})$ is close to the lower bound, if a large part of the cut locus $CL(\alpha_{j})$ of $\alpha_j$ is contained in the corresponding Q-piece $\mathcal{Q}_i$ and if $\ell(\alpha_j)\cdot |{\bf tw}_j|$ is small. The first condition is fulfilled if the length $\ell(\beta_i)$ of the boundary geodesic $\beta_i$ is small, the second if both $\ell(\alpha_j)$ and $|{\bf tw}_j|$ are small. This justifies the choice of $\alpha_i$ in inequality (\ref{eq:parlier1}). It is noteworthy that for $|{\bf tw}_j|=0$ and $\ell(\beta_i)$ small the estimates are almost sharp, independent of the length $\ell(\alpha_j)$.\\
The estimates for the entries of $P_S$ are linear combinations of the upper and lower bounds for the energies of dual harmonic forms. Hence these estimates are good if all Fenchel-Nielsen coordinates involved are small. Note that by \cite{bse2} there exists a canonical basis for a Riemann surface of genus $g$, where the largest element is of order $g$. Hence, at least the condition on the length of the geodesics involved can in principle be satisfied for small $g$.\\
The advantage of the method is that information about the geometry of the surface can be incorporated. Suppose, for example, that the geometry of $\mathcal{Y}_1$, the surface of signature $(0,3)$, or \textit{Y-piece}, attached to the Q-piece $\mathcal{Q}_1$ is known. Then for $j \in \{1,2,12\}$ the cut locus $CL(\alpha_j) \cap (\mathcal{Q}_1 \cup \mathcal{Y}_1)$ can be calculated. Incorporating this information, we obtain better estimates for the corresponding entries of the period Gram matrix. Information about isometries of the surface can also be incorporated. This is shown in \textbf{Example 3.2}.

\section{Preliminaries}

Many calculations presented in the following sections rely on the embedding of topological tubes around simple closed geodesics of Riemann surfaces into hyperbolic cylinders and subsequent approximations and calculations in Fermi coordinates. These concepts are presented in Section 2.1. Then in Section 2.2 the definition of the Fenchel-Nielsen coordinates used throughout the paper is given.

\subsection{Fermi coordinates and capacity estimates}

The Poincar\'{e} model of the hyperbolic plane is the following subset of the complex plane $\C$,
\[
   \Hy = \{ z = x + i y \in \C \mid y > 0 \}
\]
with the hyperbolic metric $ds^2 = \frac{1}{y^2}(dx^2 + dy^2)$.\\
\textit{Fermi coordinates  $\psi$, with baseline $\eta$ and base point $q_0$}, are defined as follows: the Fermi coordinates are a bijective parametrization of $\Hy$
\[
  \psi : \R^2 \rightarrow \Hy, \psi: (t,s) \mapsto \psi(t,s),
\]
where  $\psi(0,0)=q_0$. Each point $q=\psi(t,s) \in \Hy$ can be reached by starting from the base point $q_0$ and moving along $\eta$, the directed distance $t$ to $\psi(t,0)$. There is a unique geodesic, $\nu$, intersecting $\eta$ perpendicularly in $\psi(t,0)$. From $\psi(t,0)$ we now move along $\nu$ the directed distance $s$ to $\psi(t,s)$.\\
A \textit{hyperbolic cylinder} $C$ or shortly \textit{cylinder} is a set isometric to
\[
\{ \psi(t,s) \mid (t,s) \in [0,b]\times [a_1,a_2]\} \mod \{ \psi(0,s) = \psi(b,s) \mid s \in  [a_1,a_2]\},
\]
with the induced metric from $\Hy$. The baseline of $C$ is the simple closed geodesic $\gamma$ in $C$, which has length $\ell(\gamma)=b$.\\
Consider a cylinder $C$. Let $U \subset C$ be a set and $F \in \lip(\bar{U})$ a Lipschitz function on the closure of $U$. Let $G$ be the metric tensor with respect to the Fermi coordinates. Then the \textit{energy $E_U(F)$ of $F$ on $U$} is given by
\[
    E_U(F) = \iint\limits_{\psi^{-1}(U) } {\left\| {D(F \circ \psi)} \right\|_{G^{-1}}^2 \sqrt{\det(G)}}.
\]
Using Fermi coordinates, we obtain $E_U(F)$ with $F \circ \psi =f$:
\begin{equation}
E_U(F)=\iint\limits_{\psi^{-1}(U)}{\frac{1}{\cosh(s)}\frac{\partial f(t,s)}{\partial t}^2+\cosh(s)\frac{\partial f(t,s)}{\partial s}^2 \,ds \, dt } \geq \iint\limits_{\psi^{-1}(U)}\cosh(s)\frac{\partial f(t,s)}{\partial s}^2.
\label{eq:energy}
\end{equation}
The \textit{capacity} $\capa(R)$ of an annulus $R \subset C$ is given by
\[
   \capa(R)= \mathop {\inf }\{E_R(F) \mid \{ F \in \lip(\bar{R}) \mid F |_{\partial_1 R} = 0 , F |_{\partial_2 R} = 1 \} \}.
\]
In \cite{mu1}, we obtain general upper and lower bounds on the capacity of annuli on a cylinder of constant curvature. These annuli are obtained by a continuous deformation of the cylinder itself. A lower bound is obtained by determining explicitly the function that satisfies the boundary conditions of the capacity problem on $R$ and minimizes the last integral in the above inequality (\ref{eq:energy}). If the annulus $R \subset C$ is given in Fermi coordinates by
\[
R=\psi\{ (t,s) \mid s \in [a_1(t),a_2(t)],t \in [0,\ell(\gamma)] \},
\]
where $a_1(\cdot)$ and $a_2(\cdot)$ are piecewise differentiable functions with respect to $t$. Then by \cite{mu1}, \textbf{Theorem 4.1} we have:

\begin{thm}
There exists a test function $F^{test} \in \lip(R)$,  such that for $H(s)= 2\arctan(\exp(s))$ and $q_i(t)=\frac{\partial H(s_0)}{\partial s}|_{s_0=a_i(t)}\cdot a_i'(t)$ for $i \in \{1,2\}$, the capacity of $R$ satisfies:
\[
\int\limits_{t=0}^{\ell(\gamma)} {\frac{1+\frac{q_1(t)^2+q_1(t)q_2(t)+q_2(t)^2}{3}}{H(a_2(t))-H(a_1(t))} \, dt }   = E(F^{test}) \geq \capa(R) \geq \int\limits_{t=0}^{\ell(\gamma)} {\frac{1}{H(a_2(t))-H(a_1(t))} \, dt}.
\]
\label{thm:capa_S2}
\end{thm}
The estimate is sharp if the boundary is constant. In this case, $a_1(t)=a_1$, $a_2(t)=a_2$, and $R$ is itself an embedded cylinder. Especially, for $a_1=-a_2$, this simplifies to 
\begin{equation}
\capa(R)=  \frac{\ell(\gamma)}{\pi - 2\arcsin((\cosh(a_1))^{-1}) }.
\label{eq:mk1}
\end{equation}
The estimate worsens, if the variation of the boundary, $\displaystyle{ \int\limits_{t = 0}^{\ell(\gamma)} {|a_1'(t)|^2 + |a_2'(t)|^2 } \,dt }$ increases.

\subsection{Y-pieces and Fenchel-Nielsen coordinates}

\begin{figure}[h!]
\SetLabels
\L(.49*.11) $\gamma_1$\\
\L(.39*.84) $\gamma_2$\\
\L(.59*.84) $\gamma_3$\\
\L(.39*.45) $c_{12}$\\
\L(.48*.70) $c_{23}$\\
\L(.58*.45) $c_{31}$\\
\endSetLabels
\AffixLabels{%
\centerline{%
\includegraphics[height=5cm,width=6cm]{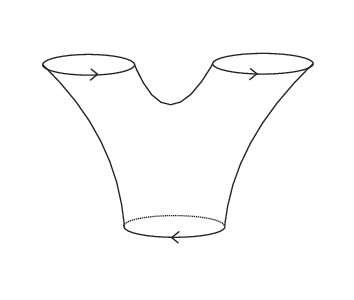}}}
\caption{A marked Y-piece $\mathcal{Y}[\gamma_1,\gamma_2,\gamma_3]$}
\label{fig:marked_Y}
\end{figure}
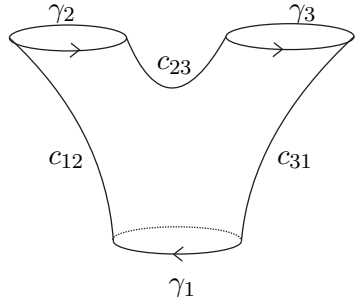

An important class of hyperbolic Riemann surfaces with geodesic boundary are the surfaces of signature $(0,3)$, or \textit{Y-pieces}. Any Riemann surface of signature $(g,n)$ can be decomposed into or built from these basic building blocks. The geometry of a hyperbolic Y-piece is determined by the length of its three boundary geodesics.\\
If $\mathcal{Y}$ is a Y-piece with boundary geodesics $\gamma_1,\gamma_2,\gamma_3$, then we can introduce a marking on $\mathcal{Y}$. The marking entails labelling the boundary components to obtain the \textit{marked Y-piece }$\mathcal{Y}[\gamma_1,\gamma_2,\gamma_3]$. For $\mathcal{Y}[\gamma_1,\gamma_2,\gamma_3]$, we introduce a standard parametrization of the boundaries as explained below.\\
Let $c_{ij}$ be the geodesic arc going from $\gamma_i$ to $\gamma_j$ that meets these boundaries perpendicularly. We set $\Sp^1 = \R \mod (t \mapsto t+1)$ and parametrize all boundary geodesics
\[
\gamma_i: \Sp^1  \rightarrow \mathcal{Y}[\gamma_1,\gamma_2,\gamma_3],  \gamma_i: t \mapsto \gamma_i(t),
\]
such that each geodesic is traversed once and with the same orientation. We parametrize the geodesics, such that  $\gamma_1(0)$ is the endpoint of $c_{31}$, $\gamma_2(0)$ is the endpoint of $c_{12}$, and $\gamma_3(0)$ is the endpoint of $c_{23}$ (see \textit{Figure \ref{fig:marked_Y}}).

Two marked Y-pieces $\mathcal{Y}$ and $\mathcal{Y}'$, that have a boundary geodesic of the same length, can be pasted together. If $\gamma_1 \subset \mathcal{Y}$  and $\gamma_1' \subset \mathcal{Y}'$ are the geodesics of equal length, then we can glue $\mathcal{Y}$ and $\mathcal{Y}'$ using the identification
\[
  \gamma_1(t)= \gamma_1'(-t +{\bf tw}) , t \in \Sp^1,
\]
where ${ \bf tw} \in \R$ is an additional constant, called the \textit{twist parameter}. We  obtain the surface
\[
     \mathcal{Y}+\mathcal{Y}' \mod (\gamma_1(t)= \gamma_1'(-t + {\bf tw}), t \in \Sp^1).
\]
If $\gamma$ is the simple closed geodesic in $\mathcal{Y}+\mathcal{Y}'$, which corresponds to $\gamma_1$ in $\mathcal{Y}$, then we call ${\bf tw}$ the twist parameter at $\gamma$.\\
Every Riemann surface $S$ of signature $(g,n)$ can be built from $2g-2 + n$ Y-pieces. The pasting scheme can be encoded in a graph $G(S)$ (see \cite{bu}, pp. 27 -30). Let $L(S)$ be the set of $3g-3 +n$ lengths of simple closed geodesics in the surface $S$, corresponding to the boundary geodesics of the Y-pieces from the construction. Let $B(S)$ be the set of $3g-3+n$ twist parameters that define the gluing of these geodesics. Then any Riemann surface $S$ can be constructed from the information provided in the triplet $(G(S),L(S),B(S))$.
\begin{defi}
$(L(S),B(S))$ is the sequence of Fenchel-Nielsen coordinates of the Riemann surface $S$.
\end{defi}
We finally note that, up to isometry, any Riemann surface can be constructed taking all twist parameters in the interval $(-\frac{1}{2},\frac{1}{2}]$ (see \cite{bu}, \textbf{Theorem 6.6.3}) and we will make use of this fact in Section 4. 

\section{Theoretical estimates for the period Gram matrix}

Let $S$ be a Riemann surface of genus $g \geq 2$, ${\rm A}$ a canonical basis, and $P_S$ the period Gram matrix of $S$ with respect to ${\rm A}$,
\[
    P_S =  \left(p_{ij}\right)_{i,j= 1,...,2g} = \left( {\int\limits_S {\sigma_i  \wedge {}^ * } \sigma_j } \right)_{i,j= 1,...,2g}.
\]
Here we first show how to obtain bounds on the diagonal entries of $P_S$ using the geometry of embedded cylinders around the elements of the canonical basis. This approach can be elaborated to obtain estimates for all entries of $P_S$. It relies on the premise that the cut locus of a given simple closed geodesic on a Riemann surface can be (at least partially) calculated.

\subsection{Estimates for the diagonal entries of $P_S$}

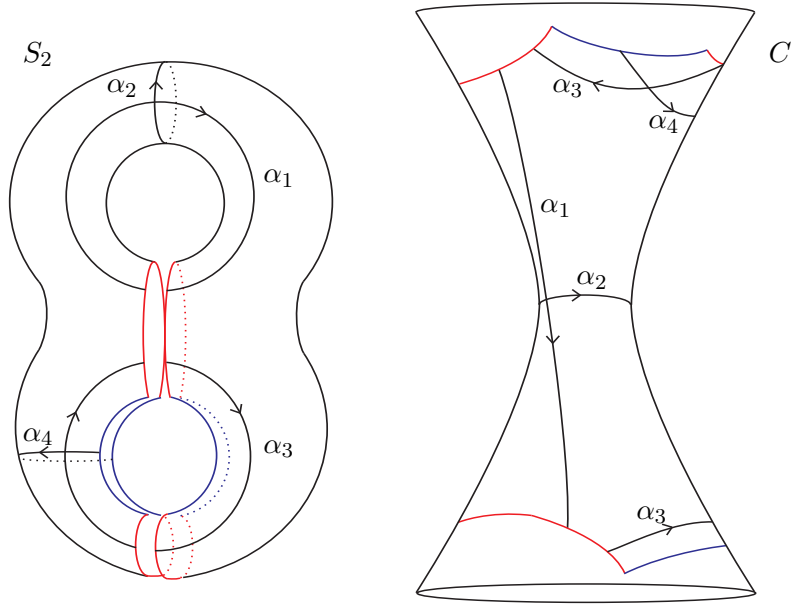
\begin{figure}[h!]
\SetLabels
\L(.20*.90) $S_{2}$\\
\L(.27*.85) $\at$\\
\L(.40*.70) $\ao$\\
\L(.40*.25) $\alpha_3$\\
\L(.20*.27) $\alpha_4$\\
\L(.82*.90) $C$\\
\L(.66*.53) $\at$\\
\L(.63*.65) $\ao$\\
\L(.71*.14) $\alpha_3$\\
\L(.64*.85) $\alpha_3$\\
\L(.72*.79) $\alpha_4$\\
\endSetLabels
\AffixLabels{%
\centerline{%
\includegraphics[height=8cm,width=10cm]{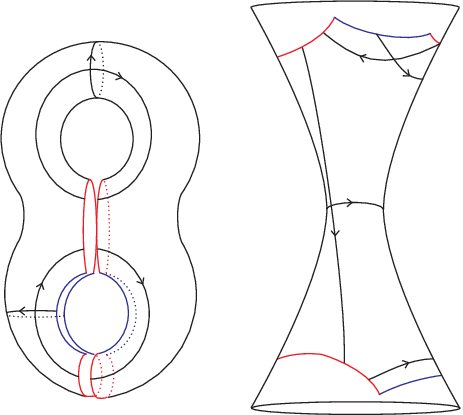}}}
\caption{Embedding of $S_2=S_{\alpha_2}$ in a cylinder around $\alpha_2$}
\label{fig:red_blue}
\end{figure}

Let $T(\alpha_{\tau(i)}) \subset S$ be a topological tube that contains the geodesic $\alpha_{\tau(i)}$ in its interior. We will see in Section 3.1.1 that the capacity of such a tube gives an upper bound for the energy of $\sigma_i$, $E(\sigma_i)=p_{ii}$. Consider without loss of generality $E(\sigma_1)=p_{11}$.\\
We will use the tube obtained by cutting open $S$ along the cut locus $CL(\alpha_2)$ of $\alpha_2$. The \textit{cut locus} of a subset $X \subset S$, $CL(X)$ is defined as follows:
\begin{equation}
   CL(X) := \{y \in S \mid \exists \gamma_{x,y},\gamma_{x',y}, \gamma_{x,y} \neq \gamma_{x',y}, \text{ with } x,x' \in X \text{ and } \dist(x,y)=\ell(\gamma_{x,y}) = \ell(\gamma_{x',y}) \},
\label{eq:cut_locus}
\end{equation}
where $\gamma_{a,b}$ denotes a geodesic arc connecting the points $a$ and $b$.
We denote by $S_X$ the surface, which we obtain by cutting open $S$ along $CL(X)$. For a set $X \subset S$, set
\begin{equation}
Z_r(X) = \{ x \in S \mid \dist(x,X) \leq r \}.
\label{eq:zrx}
\end{equation}
If $U$ is a union of disjoint simple closed geodesics $(\gamma_i)_{i=1,...,n}$, then for a sufficiently small $r$, $Z_r(U)$ consists of disjoint cylinders around these geodesics. We obtain $CL(U)$ by letting $r$ grow continuously until $Z_r(U)$ self-intersects. We stop the expansion at the points of intersection, but continue expanding the rest of the set, until the process halts. The points of intersection then form $CL(U)$. It follows from this process that the surface $S_U$, that we obtain by cutting open $S$ along $CL(U)$, can be retracted onto the union of small cylinders around $(\gamma_i)_{i=1,...,n}$. If $U=\gamma$, then $S_U$ can be embedded into an sufficiently large cylinder $C$ around $\gamma$. For more information about the cut locus, see \cite{ba}.\\
Consider an embedding of $S_{\at}=S_2$ in a cylinder $C$ (see \textit{Figure \ref{fig:red_blue}}), which, by abuse of notation, we also call $S_2$. The boundary $\partial S_2$ of $S_2 \subset C$ consists of the two connected components $\partial_1 S_2$ and $\partial_2 S_2$, which are piecewise geodesic (see \cite{ba}).
Fixing a base point $x \in S_2 \subset C$, we can construct a primitive $F_1$ of $\sigma_1$ by integrating $\sigma_1$ along paths starting from the base point $x$. As $\int\limits_{\alpha_2} \sigma_1 =0$, the value of the integral is independent of the chosen path in $S_2$. Hence, there exists a primitive $F_1$ of $\sigma_1$ on $S_2 \subset C$. Furthermore, $F_1$ is a real harmonic function, as $\sigma_1$ is a real harmonic 1-form. We recall that the value of the integral of $\sigma_1$ over a closed curve depends only on the homology class of the curve. In particular, the value of the integral is the same for two curves in the same free homotopy class. The boundary $\partial S_2 \subset C$ has two connected components, $\partial_1 S_2$ and $\partial_2 S_2$ that lie on opposite sides of $\alpha_2$. The conditions on the canonical basis ${\rm A}$ imply the following boundary conditions for $F_1$. For each point $p_1$ on $\partial_1 S_2 \subset C$, there exists a point $p_2$, such that $p_1$ and $p_2$ map to the same point $p$ on $S$, and
\[
F_1(p_2)-F_1(p_1)=0 \text{ \ \ or \ \ } F_1(p_2)-F_1(p_1)= 1.
\]
We color $p_1$ and $p_2$ blue if $F_1(p_2)-F_1(p_1)=0$ and red if  $F_1(p_2)-F_1(p_1)= 1$. If $p_2$ is on $\partial_2 S_2$ and $F(p_1) - F(p_2) = 0$ for a point $p_1$, such that $p_2 = p_1$ on $S$, then we also color $p_1$ and $p_2$ blue. We call such a decomposition a \textit{red-blue decomposition} of the cut locus (see \textit{Figure \ref{fig:red_blue}}).\\
Let $CL^{blue}(\alpha_2)$ and $CL^{red}(\alpha_2)$ denote the blue and the red parts of $CL(\at)$, both in $S$ and $S_2$. Then
\[
    CL(\at) = CL^{blue}(\alpha_2) \cup CL^{red}(\alpha_2).
\]
For the red-blue decomposition that is obtained via the cut locus $CL(\at)$, the following holds. If $p_1$ and $p_2$ are blue, then $p_1$ and $p_2$ lie both on either $\partial_1 S_2$ or $\partial_2 S_2$. If $p_1$ and $p_2$ are red, then they lie on opposite sides. This follows from the relationship of the canonical 1-forms with the intersection number of curves (see \cite{fk}, chapter III).  At the intersection of the red and the blue parts of a boundary, there exist a finite number points that are both red and blue.\\
We now connect the endpoints of two corresponding opposite red boundary segments in the red-blue decomposition of $S_2 \subset C$ with differentiable curves. Then the curves, together with the boundary segments of $S_2$, enclose a subset of $S_2$. Note that some of these curves may cross. Therefore we choose from $S_2$ a subset of curves, such that these do not mutually intersect and denote by $S_2^{red}$ the union of all enclosed areas obtained this way.

\subsubsection{Upper bound} \label{sec:upper_diag}

Let $T(\alpha_2) \subset S$ be a topological tube (with piecewise differential boundary) that contains the geodesic $\alpha_2$ in its interior. Let $\tilde{\sigma}_1$ be a closed 1-form that satisfies

\begin{equation}
\int\limits_{[\alpha_{k}]} {\tilde{\sigma}_{1} } = \delta_{1k}     \text{ \ for all \ } k \in \{1,...,2g\}.
\label{eq:uprime}
\end{equation}

Then $\sigma_1$ is the unique energy-minimizing closed 1-form satisfying the above equation. Hence, $E(\sigma_1) \leq E(\tilde{\sigma}_1)$.\\   
Let $F$ be a function, that solves the capacity problem for $T(\alpha_2)$, i.e.  
\[
F |_{\partial_1 T(\alpha_2)}=0 \text{ \ \ , \ \ } F  |_{\partial_2 T(\alpha_2)}=1 \text{ \ and \ } \capa(T(\alpha_2))=E(F). 
\]
This function is harmonic in the interior of $T(\alpha_2)$. We obtain a 1-form $\sigma'_1$ setting
\[
 \sigma'_1=
\left\{ {\begin{array}{*{20}c}
   {0}  \\
   {DF }  \\
\end{array}} \right. \text{  on  } \begin{array}{*{20}c}
   { S \backslash \{T(\alpha_2)\}  }  \\
   {T(\alpha_2)}  \\
\end{array}.
\]
By Stoke's theorem it satisfies Equation (\ref{eq:uprime}). Now, $\sigma'_1$ might not be a closed 1-form. However, there is always a closed 1-form arbitrarily near $\sigma'_1$ that also satisfies (\ref{eq:uprime}) and we may assume that $\sigma_1$ is closed. This implies that
\[
      E(\sigma_1) \leq E(\sigma'_1) = \capa(T(\alpha_2)).
\]
The above inequality is also used in \cite{bs}, where $T(\alpha_2)$ is an embedded cylinder. Setting $T(\alpha_2)=S_2$, we obtain that the capacity $\capa(S_2)$ of $S_2 \subset C$ provides an upper bound on the energy of $\sigma'_1$.\\
We obtain an upper bound on $\capa(S_2)$ by evaluating the energy of any test function $F^{test}_{1}$, that is a Lip\-schitz function on $S_2$ and satisfies the boundary conditions of the capacity problem (see \cite{tr}). As $S_2 \subset C$ is an annulus that satisfies the conditions of \textbf{Theorem \ref{thm:capa_S2}}, such a function $F^{test}_1$ is provided there:
\begin{equation}
    E(F^{test}_1) \geq \capa(S_2) \geq E(\sigma_1) = p_{11}.
\label{eq:F_1t}    
\end{equation}

\subsubsection{Lower bound} \label{sec:lower_diag}

We obtain a lower bound on $p_{11}=E(F_1)$ as explained next. Consider the set $S_2^{red}$. Remember that in each connected subset of $S_2^{red}$ there are boundary points $p_1$ and $p_2$ on opposite sides, such that $F_1(p_2)-F_1(p_1)= 1$. We obtain:
\[
   E(\sigma_1) = E(F_1) \geq \int\limits_{S_2^{red}} \| DF_1 \|_2^2.
\]
Let $I$ be a disjoint union of intervals in $\R$ and
\[
\varphi : I \times [b_1,b_2] \rightarrow S_2^{red}, \varphi : (t,s) \mapsto \varphi(t,s)
\]
a bijective function that parametrizes $S_2^{red}$ as follows:
\[
    \varphi(I\times \{b_1\})=S_2^{red} \cap \partial_1 S_2 \text{ \ \ and \ \ }   \varphi(I\times \{b_2\})=S_2^{red} \cap \partial_2 S_2
\]
and for a fixed $t_0 \in I$, $\varphi(\{t_0\}\times [b_1,b_2])$ is a differentiable curve in $S_2^{red}$, such that
\[
     F_1(\varphi(t_0,b_2))-F_1(\varphi(t_0,b_1))=1.
\]
Denote by $\mathcal{F}_1$ the set of functions
\[
\mathcal{F}_1 =\{ f: S_2^{red} \rightarrow \R \mid f \in \lip(S_2^{red}) \text{ \ and \ }  f(\varphi(t_0,b_2))-f(\varphi(t_0,b_1))=1 \text{ \ }  \forall \text{ \ }   t_0 \in I \}.
\]
We can obtain a lower bound on $p_{11}=E(\sigma_1)=E(F_1)$ if we find a function $\tilde{f}_1$, such that

\begin{equation}
\int\limits_{S_2^{red}} \|D\tilde{f}_1 \|_2^2= \min  \limits_{f \in \mathcal{F}_1} \int\limits_{S_2^{red}} \| Df \|_2^2.
\label{eq:noprojection}
\end{equation}
We call this problem the \textit{free boundary problem} for $S_2^{red}$.\\
Though this problem is quite interesting in its own right, we could not find an explicit solution.
To obtain an explicit result, we construct another lower bound based on projection of tangent vectors on curves. For an $x=\varphi(t_0,s_0) \in S_2^{red}$ denote by 
\[
pr_{{\varphi},x} : T_x(S_2^{red}) \rightarrow \{\lambda \cdot \frac{\partial \varphi(t_0,s_0)}{\partial s} \mid \lambda \in \R \}
\]
 the orthogonal projection of a tangent vector in $x$ onto the subspace spanned by $\frac{\partial \varphi(t_0,s_0)}{\partial s}$. With the help of this projection $pr_{{\varphi}} : T(S_2^{red}) \rightarrow T(S_2^{red})$ we get:
\begin{equation}
    E(F_1) \geq \int\limits_{S_2^{red}} \| DF_1 \|_2^2 \geq \int\limits_{S_2^{red}} \| pr_{\varphi}(DF_1) \|_2^2 \geq  \min \limits_{f \in \mathcal{F}_1} \int\limits_{S_2^{red}} \| pr_{\varphi}(Df) \|_2^2 = \int\limits_{S_2^{red}} \| pr_{\varphi}(Df_1) \|_2^2.
\label{eq:projection}
\end{equation}
Here, $f_1$ is a function that realizes the minimum. We have  $\int\limits_{S_2^{red}} \| DF_1 \|_2^2 = \int\limits_{S_2^{red}} \| pr_{\varphi}(DF_1) \|_2^2$, if and only if in every point $\varphi(t_0,s_0)=x \in S_2^{red}, \varphi(t_0,\cdot)$ is orthogonal to the level set of $F_1$ passing through $x$. Note that the problem of finding the function $f_1$ is in general easier than finding the function $F_1$ or $\tilde{f}_1$. We will apply these ideas to Q-pieces in Section 4. Summarizing the inequalities (\ref{eq:F_1t})-(\ref{eq:projection}) we obtain the following estimates for a diagonal entry of the period Gram matrix $p_{11}= E(\sigma_1)$:
\begin{equation}
 E(F^{test}_{1}) \geq \capa(S_{\at}) \geq  E(\sigma_1) = E(F_1) \geq \min \limits_{f \in \mathcal{F}_1} \int\limits_{S_2^{red}} \|Df \|_2^2  \geq \min \limits_{f \in \mathcal{F}_1} \int\limits_{S_2^{red}} \| pr_{\varphi}(Df) \|_2^2.
\label{eq:diag_all}
\end{equation}
Note that the upper bound differs from the lower bound. One reason for this difference is that the test function whose energy provides our upper bound has positive energy on $S_2 \backslash S_2^{red}$, whereas the energy is zero in the estimate providing the lower bound.  Another difference is due to the use of the projection along lines in the construction of the lower bound. In Section 4, we will apply these methods to a decomposition of the Riemann surface, where the elements of the canonical basis are contained in Q-pieces. There, we will see these two effects explicitly. Generalizing (\ref{eq:diag_all})  we obtain for a diagonal entry of the period matrix:
\begin{thm}\textit{(diagonal entries of $P_S$)} Let $S$ be a Riemann surface of genus $g \geq 2$ and ${\rm A} = (\alpha_i)_{i=1,...,2g}$ be a canonical homology basis and $(\sigma_i)_{i=1,...,2g}$ be the corresponding dual basis of harmonic 1-forms. Let $ S_{\alpha_k} = S_k$ be the surface obtained by cutting open $S$ along the cut locus $CL(\alpha_k)$ of $\alpha_k$ and $P_{S}= (p_{ij})_{i,j}$ be the period Gram matrix with respect to ${\rm A}$. Then
\begin{equation*}
 E(F^{test}_{i}) \geq \capa(S_{\tau(i)}) \geq  E(\sigma_i) = p_{ii} = E(F_i)  \geq \min \limits_{f \in \mathcal{F}_i} \int\limits_{S_{\tau(i)}^{red}} \|Df \|_2^2  \geq \min \limits_{f \in \mathcal{F}_i} \int\limits_{S_{\tau(i)}^{red}} \| pr_{\varphi}(Df) \|_2^2.
\end{equation*}
\label{thm:diag}
\end{thm}

\subsection{Estimates for the non-diagonal entries of $P_S$}

We now show, how we can estimate the remaining entries of the period Gram matrix $P_S$. Since $\int\limits_S {\cdot \wedge {}^ *  \cdot }$ is a scalar product, for $i \neq j$ we have by the polarization identity:
\begin{eqnarray}
\label{eq:sca_norm}
|p_{ij}|&\leq& \frac{1}{2}\left(E(\sigma_i) + E(\sigma_j) \right), \\
\label{eq:sigma_sc+}
p_{ij}&=& \frac{1}{2}\left( E(\sigma_i+\sigma_j) - E(\sigma_i) - E(\sigma_j) \right) , \text{ \ and  \ }   \\
\label{eq:sigma_sc-}
p_{ij}&=& \frac{1}{2}\left(E(\sigma_i) + E(\sigma_j) -E(\sigma_i-\sigma_j) \right).
\end{eqnarray}

We have shown how to find upper and lower bounds on $E(\sigma_i)$ and $E(\sigma_j)$. We obtain a direct estimate of $p_{ij}$ from inequality (\ref{eq:sca_norm}). However, to obtain a sharp estimate, both $E(\sigma_i)$ and $E(\sigma_j)$ must be small. We will show how to obtain better estimates of $p_{ij}$ from the following two equations. If we can find upper and lower bounds on either $E(\sigma_i+\sigma_j)$ or $E(\sigma_i-\sigma_j)$, we will obtain an estimate for $p_{ij}$. Now $\sigma_i+\sigma_j$ and $\sigma_i-\sigma_j$ satisfy the following equations on the cycles:
\begin{equation}
\label{eq:sigma_ij+}
   \int\limits_{[\alpha_k]} {\sigma_i +\sigma_j} = \delta_{ik} + \delta_{jk}  \text{ \ \ \ and \ \ }
   \int\limits_{[\alpha_k]} {\sigma_i -\sigma_j} = \delta_{ik} - \delta_{jk}  \text{ \ \ \ for all \ \ } k \in \{1,...,2g\}.
\end{equation}
There is a geodesic $\alpha$ in the free homotopy class of either $\alpha_{\tau(i)} \cdot \alpha_{\tau(j)}$ or $\alpha_{\tau(i)} (\alpha_{\tau(j)})^{-1}$ which is a simple closed curve. Applying a base change of the canonical basis, we can incorporate $\alpha$ into a new basis. This can be done, such that one of the two 1-forms $\sigma_i +\sigma_j$ and $\sigma_i -\sigma_j$ becomes an element of the new dual basis. Hence we can obtain upper and lower bounds for the energy of one of these harmonic forms using the methods from the previous subsection.\\
Since it can be difficult to explicitly parametrize a suitable geodesic $\alpha$, we will present this approach only for the case $\alpha_j=\alpha_{\tau(i)}$. We present these estimates in Section \ref{sec:jistau}. If $\alpha_j \neq \alpha_{\tau(i)}$, we will present an alternative approach in Section \ref{sec:jnottau}. We will make use of both methods in Section 4.

\subsubsection{Estimates for a non-diagonal entry $p_{i\tau(i)}$} \label{sec:jistau}

Consider without loss of generality $p_{12}$. Let $\alpha_{12}$ be the simple closed geodesic in the free homotopy class of $\ao{\at}^{-1}$. We apply the base change
\[
   {\rm A} =(\alpha_1, \alpha_{2},...,\alpha_{2g}) \rightarrow  (\alpha_{12}, \alpha_{2},...,\alpha_{2g}) = {\rm A'}.
\]
This way we obtain the dual basis $\left( {\sigma' _k } \right)_{k = 1,...,2g}$ for ${\rm A'}$, where
\[
   (\sigma_1, \sigma_1 + \sigma_{2},\sigma_3,...,\sigma_{2g}) =  (\sigma'_{1}, \sigma'_{2},\sigma'_3,...,\sigma'_{2g}).
\]
Let $F_{12}=F'_2$ be a primitive of $\sigma_1+\sigma_2=\sigma'_2$  on $S_{\alpha_{12}}=S_{12}$. We embed $S_{12}$ into a cylinder $C$ and denote this surface also by $S_{12}$. Proceeding as in the previous subsection, we obtain upper and lower bounds on $E(\sigma_1+\sigma_2)=E(\sigma'_2)$ from the geometry of $S_{12}$: 
\[
   E(F^{test}_{12}) \geq \capa(S_{12}) \geq E(\sigma_1+\sigma_2) \geq E_{S_{12}^{red}}(pr_{\varphi}(Df_{12})).
\]
Here $F^{test}_{12}$ is the test function provided by \textbf{Theorem \ref{thm:capa_S2}}, whose energy provides an upper bound for $\capa(S_{12})$ and $f_{12}$ is the function constructed analogously to $f_1$ (see inequality (\ref{eq:projection})). Substituting the estimates of $E(\sigma_{1}+\sigma_2)$, $E(\sigma_1)$, and $E(\sigma_2)$ in Equation (\ref{eq:sigma_sc+}), we obtain an upper and lower bound on $p_{12}$. These equations are summarized in \textbf{Theorem \ref{thm:ndiag}} in the following subsection.

\subsubsection{Estimates for a non-diagonal entry $p_{ij}$, where $j \neq \tau(i)$} \label{sec:jnottau}

In this case $\alpha_i$ and $\alpha_j$ do not intersect. Consider without loss of generality $p_{13}$.  Then $\alpha_{\tau(1)}= \alpha_2$ and $\alpha_{\tau(3)}= \alpha_4$. Recall that $S_{\alpha_2 \cup \alpha_4}:=S_{24}$ is the surface, which we obtain by cutting open $S$ along $CL(\alpha_2 \cup \alpha_4)$ (see (\ref{eq:cut_locus})). 

\begin{figure}[h!]
\SetLabels
\L(.43*.84) $S_{24}$\\
\L(.20*.84) $S_{24}^1$\\
\L(.30*.82) $\at$\\
\L(.40*.70) $\ao$\\
\L(.20*.41) $S_{24}^3$\\
\L(.40*.28) $\alpha_3$\\
\L(.31*.12) $\alpha_4$\\
\L(.76*.90) $C_1$\\
\L(.63*.76) $\at$\\
\L(.64*.65) $\ao$\\
\L(.76*.41) $C_3$\\
\L(.64*.15) $\alpha_3$\\
\L(.63*.28) $\alpha_4$\\
\endSetLabels
\AffixLabels{%
\centerline{%
\includegraphics[height=8cm,width=10cm]{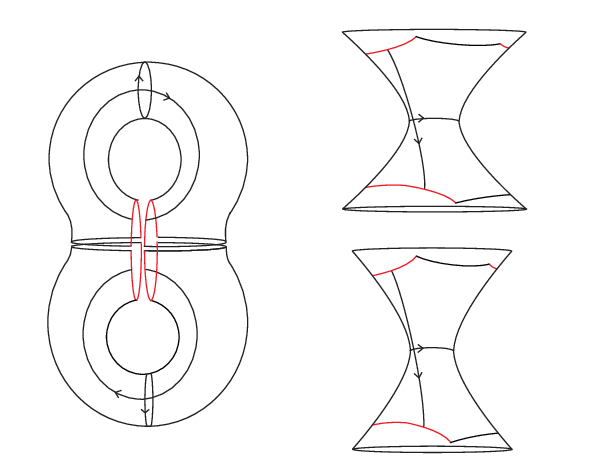}}}
\caption{Embedding of $S_{24}^i$ in a cylinder $C_i$ around $\alpha_{\tau(i)}$ for $i \in \{1,3\}$}
\label{fig:red_blue_case}
\end{figure}
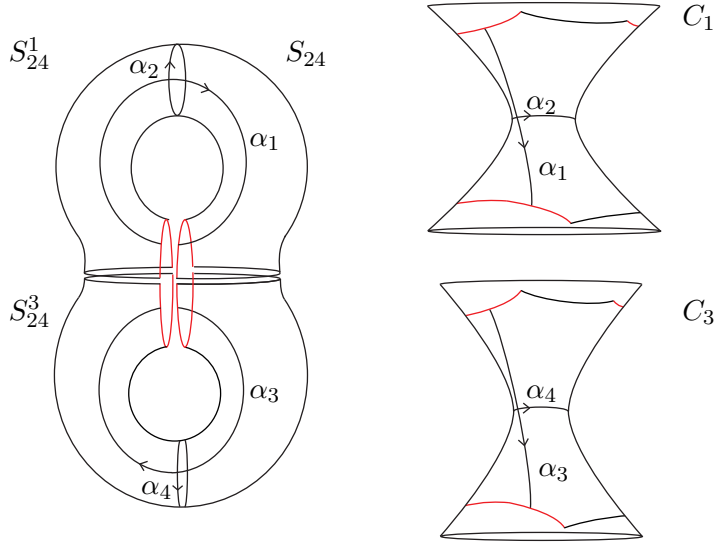

$S_{24}$ consists of two connected parts. Let $S_{24}^1 \subset S_{\at}$ be the part that contains $\alpha_2$ and let $S_{24}^3 \subset S_{\alpha_4}$ be the part that contains $\alpha_4$. We embed $S_{24}^1$ into a cylinder $C_1$ around $\alpha_2$ and $S_{24}^3$ into a cylinder $C_3$ around $\alpha_4$, and denote the embedded surfaces by the same name. Due to the relationships in Equation (\ref{eq:sigma_ij+}), $\sigma_1+\sigma_3$ has a primitive on both $S_{24}^1 \subset C_1$ and $S_{24}^3 \subset C_3$. Such a decomposition is shown in \textit{Figure \ref{fig:red_blue_case}}.

For $i \in \{1,3\}$, let $F^{test}_i$ on $S_{24}^i$ be a test function for the capacity problem on $S_{24}^i$. As in Section \ref{sec:upper_diag} we conclude
\[
        E(F^{test}_1)  \geq \capa(S_{24}^1) \geq E(\sigma_1)   \text{ \ and \ }   E(F^{test}_3)  \geq \capa(S_{24}^3) \geq E(\sigma_3).
\]      
Together, these functions naturally define a function $F^{test}_{13}$ on $S_{24}$, whose derivative $DF^{test}_{13}$ satisfies the same integral conditions on the cycles as $\sigma_1+\sigma_3$. Due to the energy-minimizing property of $\sigma_1+\sigma_3$, we obtain
\[
 E(F^{test}_{13}) = E(F^{test}_1)+ E(F^{test}_3) \geq \capa(S_{24}^1)+\capa(S_{24}^3) \geq E(\sigma_1+\sigma_3).
\]

We obtain a lower bound for $E(\sigma_{1}+\sigma_3)$  by applying the same methods used to obtain a lower bound on $E(\sigma_1)$ on $S_2$ in Section \ref{sec:lower_diag}: we obtain estimates from the red-blue decompositions induced by a primitive $F_{13}$ of $\sigma_1 + \sigma_3$ on the boundary of $S_{24}^1$ in $C_1$ and $S_{24}^3$ in $C_3$. The only difference is that we have some segments of the boundary, where the red-blue decomposition does not apply. Here we disregard these pieces in the construction of $S_{24}^{1red}$ and $S_{24}^{3red}$. As these sets are disjoint, we have to find a function $f_{13}$ that satisfies
\[
      f_{13}(p_2)-f_{13}(p_1) =F_{13}(p_2)-F_{13}(p_1)=1
\]
for all points $p_1,p_2 $ in $\partial S_{24}^{1red}$ or $\partial S_{24}^{3red}$ satisfying the above equation. Let $pr_{\varphi}$ be the projection of a vector field in the tangent space onto lines of a suitable parametrization $\varphi$ of $S_{24}^{1red}$ and $ S_{24}^{3red}$. Note that this parametrization can be extended naturally to $S_{\alpha_4}^{red}$ and $S_{\alpha_2}^{red}$. We let $f_{13}$ be a function that satisfies the above equation and minimizes the projected energy $E(pr_{\varphi}(D\cdot))$ on $S_{24}^{1red} \cup S_{24}^{3red}$. Then
\[
 E(\sigma_1+\sigma_3) \geq E_{S_{24}^{red}}(pr_{\varphi}(Df_{13})), \text{ \ where \ }   S_{24}^{red}= S_{24}^{1red} \cup  S_{24}^{3red}.
\]
Substituting the estimates for $E(\sigma_{1}+\sigma_3)$, $E(\sigma_1)$, and $E(\sigma_3)$ in Equation (\ref{eq:sigma_sc+}), we obtain upper and lower bounds on $p_{13}$:
\begin{eqnarray}
\label{eq:qilnull1}
p_{13}  &\leq& \frac{1}{2} \left( \capa(S_{24}^1)+\capa(S_{24}^3) - E_{S_{\at}^{red}}(pr_{\varphi}(Df_{1})) - E_{S_{\alpha_4}^{red}}(pr_{\varphi}(Df_{3})) \right) > 0 \text{ \ and \ } \\
 p_{13}  &\geq& \frac{1}{2} \left( E_{S_{24}^{red}}(pr_{\varphi}(Df_{13})) - \capa(S_{\at})- \capa(S_{\alpha_4}) \right) < 0.
\label{eq:qilnull}
\end{eqnarray}
In the above equation, for $i \in \{1,3\}$, $f_i$ is the minimizing function corresponding to a primitive $F_i$ of $\sigma_i$ on $S_{\alpha_{\tau(i)}}^{red}$ given in \textbf{Theorem \ref{thm:diag}}. That our estimate for $p_{13}$ in the first inequality is bigger than zero can be seen as follows. By construction, we have $S_{24}^1 \subset S_{\alpha_2}$  and $S_{24}^3 \subset S_{\alpha_4}$. Now, if an annulus $R_1$ is contained in an annulus $R_2$, then $\capa(R_1) \geq \capa(R_2)$.
Hence
\[
\capa( S_{24}^1 ) \geq \capa(S_{\alpha_2}) > E_{S_{\at}^{red}}(pr_{\varphi}(Df_{1})) \text{ \ and \ } \capa( S_{24}^3 ) \geq \capa(S_{\alpha_4}) > E_{S_{\alpha_4}^{red}}(pr_{\varphi}(Df_{3})),
\]
from which follows the last inequality in (\ref{eq:qilnull1}). \\
It follows furthermore from the boundary conditions of the functions $F_1, F_3$, and $F_{13}$ that $\partial S_{24}^{1red} \subset \partial S_{\alpha_2}^{red} \text{ \ \ and \ \ } \partial S_{24}^{3red} \subset \partial S_{\alpha_4}^{red}.$ Hence
\[
E_{S_{24}^{red}}(pr_{\varphi}(Df_{13})) = E_{S_{24}^{1red}}(pr_{\varphi}(Df_{13})) + E_{S_{24}^{3red}}(pr_{\varphi}(Df_{13})) \leq  E_{S_{\at}^{red}}(pr_{\varphi}(Df_{1})) + E_{S_{\alpha_4}^{red}}(pr_{\varphi}(Df_{3})).
\]
Now the second inequality in (\ref{eq:qilnull}) follows from this inequality and the fact that \[
E_{S_{\at}^{red}}(pr_{\varphi}(Df_{1})) < \capa(S_{\alpha_2}) \text{ \ and \ }  E_{S_{\alpha_4}^{red}}(pr_{\varphi}(Df_{3})) < \capa(S_{\alpha_4}).
\]
Using this approach, we can only obtain optimal estimates if $p_{13}$ is close to zero. This is due to the fact that we do not have full information of the boundary values on our tubes $S_{24}^{3}$ and $S_{24}^{1}$. This estimate is however better than the one obtained from Equation (\ref{eq:sca_norm}). Note that by \cite{bmm} the value of $p_{13}$ is  close to zero, if $\at$ and $\alpha_4$ are separated by a small separating simple closed geodesic $\gamma$.
Generalizing the notation used in this subsection, we summarize its results in the following theorem:
\begin{thm}(\textit{non-diagonal entries of $P_S$}) Let $S$ be a Riemann surface of genus $g \geq 2$ and ${\rm A} = (\alpha_i)_{i=1,...,2g}$ be a canonical homology basis and $(\sigma_i)_{i=1,...,2g}$ be the corresponding dual basis of harmonic 1-forms. Let $P_{S}= (p_{ij})_{i,j}$ be the period Gram matrix with respect to ${\rm A}$. Let $\alpha_{i\tau(i)}$ be a simple closed geodesic in the free homotopy class of $\alpha_i(\alpha_{\tau(i)})^{-1}$. Let $S_{\alpha_k} = S_k$ be the surface obtained by cutting open $S$ along the cut locus $CL(\alpha_k)$ of $\alpha_k$. Then we obtain for the non-dialgonal entry $p_{i \tau(i)}$ 
\begin{eqnarray*}
\frac{1}{2} \left( E_{S_{i\tau(i)}^{red}}(pr_{\varphi}(Df_{i \tau(i)})) - \capa(S_{i}) - \capa(S_{\tau(i)}) \right) \leq p_{i \tau(i)} \leq \\ \frac{1}{2} \left( \capa(S_{i \tau(i)}) - E_{S_{i}^{red}}(pr_{\varphi}(Df_{i})) -  E_{S_{\tau(i)}^{red}}(pr_{\varphi}(Df_{\tau(i)}) \right).
\end{eqnarray*}
Let $p_{ij}$ be a non-diagonal entry, where $j \neq \tau(i)$. Then
\begin{eqnarray*}
p_{ij}  &\leq& \frac{1}{2} \left( \capa(S_{\tau(i)\tau(j)}^i)+\capa(S_{\tau(i)\tau(j)}^j) - E_{S_{\tau(i)}^{red}}(pr_{\varphi}(Df_{i})) - E_{S_{\tau(j)}^{red}}(pr_{\varphi}(Df_{j})) \right) > 0 \text{ \ and \ } \\
 p_{ij}  &\geq& \frac{1}{2} \left( E_{S_{\tau(i)\tau(j)}^{red}}(pr_{\varphi}(Df_{ij})) - \capa(S_{\tau(i)})- \capa(S_{\tau(j)}) \right) < 0.
\end{eqnarray*}
\label{thm:ndiag}
\end{thm}

\subsection{Examples}

We now give two examples to demonstrate the weaknesses and strengths of our method. We first show that the energy of a dual harmonic form can be lower than the capacity of a cylinder of even infinite length. The upper bound on $p_{22}$ in the following example is due to Peter Buser.\\

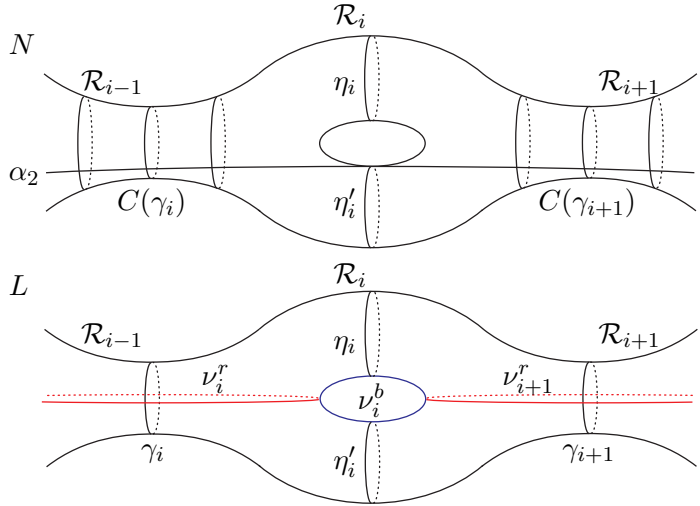
\begin{figure}[h!]
\SetLabels
\L(.20*.86) $N$\\
\L(.26*.80) $\mathcal{R}_{i-1}$\\
\L(.29*.60) $C(\gamma_i)$\\
\L(.64*.60) $C(\gamma_{i+1})$\\
\L(.47*.80) $\eta_i$\\
\L(.47*.60) $\eta'_i$\\
\L(.47*.91) $\mathcal{R}_{i}$\\
\L(.69*.80) $\mathcal{R}_{i+1}$\\
\L(.20*.46) $L$\\
\L(.20*.65) $\at$\\
\L(.26*.38) $\mathcal{R}_{i-1}$\\
\L(.31*.19) $\gamma_i$\\
\L(.66*.19) $\gamma_{i+1}$\\
\L(.47*.37) $\eta_i$\\
\L(.47*.17) $\eta'_i$\\
\L(.49*.27) $\nu^b_i$\\
\L(.36*.31) $\nu^r_i$\\
\L(.61*.31) $\nu^r_{i+1}$\\
\L(.47*.48) $\mathcal{R}_{i}$\\
\L(.69*.38) $\mathcal{R}_{i+1}$\\
\endSetLabels
\AffixLabels{%
\centerline{%
\includegraphics[height=8cm,width=10cm]{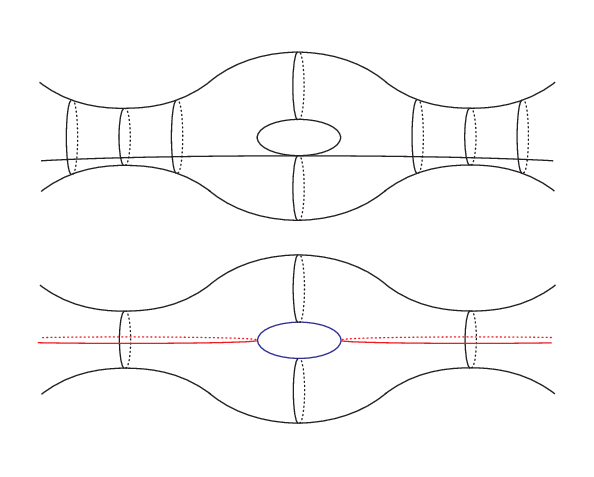}}}
\caption{Building blocks for the surfaces $N$ and $L$ of genus $g$}
\label{fig:necklace}
\end{figure}

\textbf{Example 3.1}  For comparison we briefly review the example of the necklace surface given in \cite{bse1}. Let $\mathcal{Y}$ be a Y-piece, a surface of signature $(0,3)$. Let $\gamma, \eta$ and $\eta'$ be its boundary geodesics, such that $\eta$ and $\eta'$ have equal length. We paste two copies of $\mathcal{Y}$ along $\eta$ and $\eta'$ to obtain $\mathcal{R}$ of signature $(1,2)$. As shown in \textit{Figure \ref{fig:necklace}}, the necklace surface $N$ of genus $g$ is obtained by pasting together $g-1$ copies $\mathcal{R}_1,...,\mathcal{R}_{g-1}$ of a building block $\mathcal{R}$. The free boundary of $\mathcal{R}_{g-1}$ is pasted along $\gamma_1$ of $\mathcal{R}_1$ to obtain a ring. In this example, the twist parameter for any pasting can be chosen arbitrarily.\\
By the collar lemma (see \cite{bu}, p. 106), each $\gamma_i$ has a collar of width $w_{\gamma}$, where

\begin{equation*}
   w_{\gamma} \geq \arcsinh\left(\frac{1}{\sinh(\frac{\ell(\gamma)}{2})}\right).
\end{equation*}

Let ${\rm A} = (\alpha_i)_{i=1,...,2g}$ be a canonical basis, such that $\alpha_1=\gamma_1$ and $\alpha_{\tau(1)}=\alpha_2$ is a simple closed geodesic that intersects all $(\gamma_i)_{i=1,...,g-1}$ exactly once. Let $P_S$ be the corresponding period Gram matrix. We will examine the upper bound on the entry $p_{22}=E(\sigma_2)$.\\
Following our method, we have to embed $N_1=N_{\alpha_1}$ into a cylinder $C_1$ and have to evaluate $\capa(N_1)$. Now if an annulus $R_1$ is contained in an annulus $R_2$, then $\capa(R_1) \geq \capa(R_2)$,
and hence $\capa(N_1) \geq \capa(C_1)$. From Equation (\ref{eq:mk1}), it follows that the capacity of the cylinder $C_1$ with baseline of length $\ell(\alpha_1) = \ell(\gamma)$ of infinite width is not zero. We obtain
\begin{equation}
   \capa(N_1) \geq  \capa(C_1) = \frac{\ell(\gamma)}{\pi} = \frac{\ell(\alpha_1)}{\pi}.
\label{eq:inf_cyl}
\end{equation}
We now give another estimate for the energy of $\sigma_2$ with the help of a test form $s_2$. This approach applies only to this example. To this end consider the collar $C(\gamma_i)$ of a $\gamma_i$. On each $C(\gamma_i)$  set $s_2 =DF_2$, where $F_2$ is the real harmonic function that has value $0$ on one boundary of $C(\gamma_i)$ and $\frac{1}{g-1}$ on the other. We set $s_2=0$ on $S \backslash \bigcup \limits_{i=1}^{g-1}{C(\gamma_i)}$. Then $s_2$ is arbitrarily close to a closed form that satisfies the same conditions on the elements of ${\rm A}$ as $\sigma_2$ and we have
\begin{equation}
E(\sigma_2) < E(s_2) \leq (g-1) \cdot   \frac{(g-1)^{-2}\cdot \ell(\gamma)}{\pi-2\arcsin\left({\frac{1}{\cosh(w_{\gamma})}}\right)} = \frac{c_{\alpha_1}}{g-1},
\label{eq:necklace}
\end{equation}
where $w_{\gamma}$ is bounded from below by the collar lemma.\\
The lower bound on $p_{22}$ follows from a different source. Let $\|\cdot\|_s$ be the stable norm for $H_1(S,\R)$ (see \cite{mm}, pp. 1,2 for details).
Let $J: H^1(S,\R) \rightarrow H_1(S,\R)$ be the Poincar\'e duality map. It follows from the definition in \cite{mm}, p. 4 and the corresponding \textbf{Lemma 2.2} that in \cite{mm} the map $J$ satisfies
\begin{equation}
    \int_{[\alpha]} \sigma = \ir(J(\sigma),[\alpha]) \text{ \ \ for all \ \ } [\alpha] \in H_1(S,\R), \sigma \in H^1(S,\R). 
\label{eq:poincare_dual}    
\end{equation}
Combining \textbf{Theorem 1.1} and \textbf{Lemma 2.1} of \cite{mm}, we obtain for any $[\alpha] \in H_1(S,\R)$:  
\begin{equation}
   \frac{\|[\alpha]\|^2_s}{4\pi(g-1)} \leq  E(J^{-1}([\alpha])).
\label{eq:result_MM}  
\end{equation}
It follows from the relation between the integration over cycles and the intersection form that with respect to a canonical basis 
\[ J(\sigma_2)= [\alpha_1] \text{ \ \ and \ \ }  J(\sigma_1)= -[\alpha_2]. 
\]  
This can be deduced from Equation (\ref{eq:poincare_dual}), see also \cite{jo}, Chapter 5.1 for more details. Furthermore  $\|[\alpha_1]\|_s$ is the length of a shortest multicurve in the homology class of $\alpha_1$.  We obtain from (\ref{eq:result_MM}):
\begin{equation}
     \frac{\|[\alpha_1]\|^2_s}{4\pi(g-1)} \leq E(\sigma_2).
\label{eq:p22low}  
\end{equation}
In total we obtain from Equation (\ref{eq:inf_cyl}), (\ref{eq:necklace}) and (\ref{eq:p22low}):
\begin{lem}
 Let $N_1$ be the surface obtained by cutting open the necklace surface $N$ along the cut locus $CL(\alpha_1)$ of $\alpha_1$. Then
\[
    \frac{c_{\alpha_1}}{g-1}   \geq  E(\sigma_2) \geq \frac{\|[\alpha_1]\|^2_s}{4\pi(g-1)} \text{ \ \ and \ \ } \capa(N_1) \geq E(\sigma_2), \text{ \ \ but \ \ }   \capa(N_1) \geq \frac{\ell(\alpha_1)}{\pi},
\]
where $c_{\alpha_1}$ is a factor that depends only on the fixed length $\ell(\ao)$ of $\alpha_1$ and $\|[\alpha_1]\|_s$ is the length of a shortest multicurve in the same homology class as $\alpha_1$.
\label{thm:necklace2}
\end{lem}
Hence, $E(\sigma_2)$ is of order $\frac{1}{g}$ and goes to zero as $g$ goes to infinity. Our upper bound, on the contrary, is always bigger than the constant $\frac{\ell(\gamma)}{\pi}$. This shows that there exist examples where our upper bound can not be of the right order. This might be due to the fact that the projection of $CL^{blue}(\alpha_1)$ onto $\alpha_1$ can attain almost the length of $\alpha_1$. Hence, as $CL^{blue}(\alpha_1)$ is large, $E(\sigma_2)$ might be small. \textbf{Theorem \ref{thm:necklace}} then follows from the lemma above.\\
\\
\textbf{Example 3.2} For our second example we construct a linear surface $L$ of genus $g$. This example belongs to the class of M-curves described in \cite{bsi}. In this construction, we use Y-pieces $\mathcal{Y}$, where the length of $\eta$ and $\eta'$ is large. We construct $\mathcal{R}$ from two copies of these Y-pieces as in the previous example, however, here the twist parameter in the two pastings is zero. To construct a surface $L$ of genus $g$, we paste together $g-2$ copies $\mathcal{R}_2,...,\mathcal{R}_{g-1}$ along the $\gamma_i$ (see \textit{Figure \ref{fig:necklace}}). Then, we take two copies of $\mathcal{Y}$, $\mathcal{Y}_1$, and $\mathcal{Y}_{g}$ and paste each together along $\eta$ and $\eta'$ to obtain $\mathcal{Q}_1$, and $\mathcal{Q}_{g}$, respectively. For $i\in \{1,g\}$, let $\eta_i$ denote the image of $\eta$ in $\mathcal{Q}_i$. Then we paste $\mathcal{Q}_1$ and $\mathcal{Q}_{g}$ on each side of $\mathcal{R}_2$ and $\mathcal{R}_{g-1}$, respectively. Again, the twist parameter for any pasting is zero. \\
Let ${\rm A} = (\alpha_i)_{i=1,...,2g}$ be a canonical basis, such that $\alpha_1=\eta_1$ and $\alpha_2$ is the unique simple closed geodesic in $\mathcal{Q}_1$ that intersects $\alpha_1$ perpendicularly. Let $P_L$ be the corresponding period matrix. We now show that in this case, the upper bound for $p_{22}=E(\sigma_2)$ is optimal. Therefore we use the symmetries of the surface $L$.\\
To this end, we first determine the cut locus $CL(\alpha_1)$ of $\alpha_1$. Set $\nu^b_1=\alpha_2$ and for $i \in \{2,...,g-1\}$ let $\nu^b_i \subset \mathcal{R}_i$ be the simple closed geodesic that intersect $\eta_i$ and $\eta'_i$ perpendicularly (see \textit{Figure \ref{fig:necklace}}). For $i \in \{2,...,g\}$ let $\nu^r_i$ be the simple closed geodesic that intersects the geodesic $\gamma_i$ and $\nu^b_{i-1}$  perpendicularly.  Set $\nu^r_{g+1} = \eta_g$ and let $\nu^b_g$ be the simple closed geodesic in $\mathcal{Q}_{g}$, intersecting $\eta_g$ perpendicularly.
\begin{cla}
The cut locus $CL(\alpha_1)= CL(\alpha_1)^{red} \cup  CL(\alpha_1)^{blue}$ of $\alpha_1$ in $L$ consists of the sets
\begin{equation*}
    CL(\alpha_1)^{red}= \{\nu^r_2,...,\nu^r_{g+1}\}  \text{ \ \ and \ \ } CL(\alpha_1)^{blue}= \{\nu^b_2,...,\nu^b_{g}\}.
\end{equation*}
\end{cla}
\begin{proof} To prove this claim one can use the symmetries of the surface. The proof is elementary and it is therefore left to the reader.
\end{proof}
We now show that our capacity estimate for $p_{22} = E(\sigma_2)$ is almost sharp. To this end we consider the isometries $\phi_1$, $\phi_2$ and $\phi$ in $\iso(L)$.
\begin{itemize}
\item[-]
Let $\phi_1 \in \iso(L)$ be the hyperelliptic involution that fixes $CL(\alpha_1)$ as a set, such that for all $i \in \{2,..,g-1\}$: $\,\,\, \phi_1(\eta_i)=\eta'_i$.
\item[-]
Let $\phi_2 \in \iso(L)$ be the isometry that fixes $CL(\alpha_1)$ as a set and all $\nu^b_i$ point-wise.
\item[-]
Set $\phi = \phi_1 \circ \phi_2$.
\end{itemize}
$\phi$ is the isometry that maps any point $q$ in $CL(\alpha_1)^{red} \subset S_1$ to the corresponding point $q'$ in the red-blue composition induced by $\sigma_2$. Consider a primitive $F_2$ of $\sigma_2$ on $L_{\alpha_1}=L_1$. $F_2 \circ \phi_2$ is a harmonic function, whose derivative $D(F_2 \circ \phi_2)$ defines a 1-form $\sigma'_2$ on $L$. $\sigma'_2$ satisfies the same conditions on the cycles as $\sigma_2$. Due to the uniqueness of $\sigma_2$, $\sigma'_2=\sigma_2$. In the same way $1-F_2 \circ \phi_1$ is a harmonic function, whose derivative $-D(F_2 \circ \phi_1)$ defines a 1-form $\sigma''_2$ on $L$ that satisfies the same integral conditions on the cycles as $\sigma_2$. This leads to $\sigma''_2=\sigma_2$. By choosing an appropriate additive constant, we obtain:
\begin{eqnarray*}
F_2 \circ \phi_2 &=& F_2   \text{ \ \ and \ \ } \\
1-F_2 \circ \phi_1&=&F_2 \Rightarrow  1-F_2 = F_2 \circ \phi_1.
\end{eqnarray*}
Now, for any $q$ on one side of  $CL(\alpha_1)^{red} \subset L_1$,  we have 
\[
1 = F_2(q)-F_2(\phi(q))=F_2(q)-F_2((\phi_1 \circ \phi_2)(q))).
\]
Using the two equations above this yields $F_2(q) - (1-F_2(\phi_2(q)))= 1$ or likewise $2F_2(q)=2$, hence $F_2(q)=1$. As $F_2(q)-F_2(\phi(q))=1$ it follows that $F_2(\phi(q)) = 0$. In total we obtain:
\[
    F_2(\phi(q)) = 0  \text{ \ \ and \ \ } F_2(q) = 1.
\]
Hence, the red parts of the boundary satisfy the conditions for the capacity problem. Consider the two boundary geodesics $\eta$ and $\eta'$ of our building block $\mathcal{Y}$. If $\ell(\eta) = \ell(\eta')$ is large, then it follows from hyperbolic geometry that the curves $(\nu^b_i)_{i=2,...,g}$ are arbitrarily small. The limit case is a surface $L^*$ with $2(g-1)$ cusps. The 1-form $\sigma_2$ is, however, well-defined on $L^*$ (see also \cite{bmm}). As the harmonic form $\sigma_2$ depends continuously on the domain, we obtain for small $(\nu^b_i)_{i=2,...,g}$:
\[
    p_{22} = E(\sigma_2) = \capa(L_1)- \epsilon_L,
\]
where $\epsilon_L > 0$ depends on the geometry of $L$ and may become arbitrarily small. Hence our upper bound for a diagonal entry of $P_L$ is sharp. We have shown:

\begin{lem} Let $L_1$ be the surface obtained by cutting open the linear surface $L$ along the cut locus $CL(\alpha_1)$ of $\alpha_1$. Let $(\sigma_j)_{j=1,...,2g}$ be the dual basis of harmonic forms with respect to ${\rm A}$. Then 
\[
     E(\sigma_2) = \capa(L_1) - \epsilon_L,
\]
where $\epsilon_L > 0$ depends on the geometry of $L$ and may become arbitrarily small.
\label{thm:linear2}
\end{lem}
\textbf{Theorem \ref{thm:linear}} then follows from the lemma above.

\section{Estimates for the period Gram matrix based on Q-pieces}

We note that all hyperbolic trigonometric identities in this section can be found in \cite{bu} p. 454.\\
In this section we present practical estimates for the period Gram matrix, based on the Fenchel-Nielsen coordinates of Q-pieces containing the paired curves of a canonical basis. Under this condition, the cut loci of these curves can be (at least partially) calculated.\\
More precisely, let $S$ be a Riemann surface of genus $g \geq 2$. Let $(\mathcal{Q}_i)_{i=1,3,...,2g-1} \subset S$ be a set of Q-pieces, whose interiors are pairwise disjoint. Let $\beta_i$ be the boundary geodesic of $\mathcal{Q}_i$, $\alpha_i$ an interior simple closed geodesic, and ${\bf tw}_i \in  (-\frac{1}{2},\frac{1}{2}]$ the twist parameter at $\alpha_i$. The geometry of $\mathcal{Q}_i$ is determined by the triplet $(\ell(\beta_i),\ell(\alpha_i),{\bf tw}_i)$.\\
Now fix an $i \in  \{1,3,...,2g-1\}.$ Let $\alpha_{\tau(i)} \subset \mathcal{Q}_i$ be a simple closed geodesic that intersects $\alpha_i$ exactly once, and $\alpha_{i\tau(i)} \subset \mathcal{Q}_i$ the simple closed geodesic in the free homotopy class of $\alpha_i(\alpha_{\tau(i)})^{-1}$. For $j \in \{i,\tau(i),i\tau(i)\}$, let
\begin{itemize}
\item[-] $\beta_j=\beta_i$ be the boundary geodesic
\item[-] ${\bf tw}_j$ the twist parameter at $\alpha_j$
\item[-] $FN_j := (\ell(\beta_j),\ell(\alpha_j),{\bf tw}_j)$ the corresponding Fenchel-Nielsen coordinates of $\mathcal{Q}_i$.
\end{itemize}
In \textbf{Lemma \ref{thm:FNj}} we show how to find a suitable geodesic $\alpha_{\tau(i)}$ that intersects $\alpha_i$ once and how to calculate $FN_{\tau(i)}$ and $FN_{i\tau(i)}$ from $FN_{i}$. This enables us to state estimates for all entries of the period Gram matrix $P_S$ of $S$ based on the $3g$ Fenchel-Nielsen coordinates $(FN_i)_{i=1,3,...,2g-1}$:

\begin{thm} Let $S$ be a Riemann surface of genus $g \geq 2$ and $(\mathcal{Q}_i)_{i=1,3,...,2g-1} \subset S$ be a set of Q-pieces, whose interiors do not mutually intersect. If $\mathcal{Q}_i$ is given in the Fenchel-Nielsen coordinates $FN_i=(\ell(\beta_i),\ell(\alpha_i),{\bf tw}_i)$, where $\alpha_i$ is an interior simple closed geodesic, such that $\cosh(\frac{\ell(\alpha_i)}{2}) \leq \cosh(\frac{\ell(\beta_i)}{6}) + \frac{1}{2}$.\\
Then there is a simple closed geodesic $\alpha_{\tau(i)} \subset \mathcal{Q}_i$, and a simple closed geodesic $\alpha_{i\tau(i)}$ in the free homotopy class of $\alpha_i(\alpha_{\tau(i)})^{-1}$, and the following functions
\begin{eqnarray*}
   f^u, f^l &:& \R^+ \times \R^+ \times (-\frac{1}{2},\frac{1}{2}] \rightarrow \R^+ \text{ \ (see Section 4.4) \ }\\
     f^u&:&FN_{j} \mapsto f^u(FN_{j}) \text { \ \ and \ \ }  f^l:FN_{j} \mapsto f^l(FN_j),
\end{eqnarray*}
that provide upper and lower bounds for all entries of the corresponding period Gram matrix $P_S=(p_{ij})_{i,j}$ as follows. For a diagonal entry $p_{ii}$, we have:
\[
    f^l(FN_{\tau(i)}) \leq p_{ii}   \leq  f^u(FN_{\tau(i)}).
\]
For a non-diagonal entry $p_{i\tau(i)}$, we have:
\begin{eqnarray*}
p_{i\tau(i)} \leq \frac{1}{2}\left(f^u(FN_{i\tau(i)}) - f^l(FN_{\tau(i)}) - f^l(FN_i)\right) \text{ \ and \ }  \\
p_{i\tau(i)} \geq  \frac{1}{2}\left(f^l(FN_{i\tau(i)}) -  f^u(FN_{\tau(i)}) -    f^u(FN_i) \right).
\end{eqnarray*}
For a non-diagonal entry $p_{ik}$, where $k \neq \tau(i)$, we have:
\[
0 \leq |p_{ik}| \leq \frac{1}{2}\left(  f^u(FN_{\tau(i)}) +  f^u(FN_{\tau(k)}) -  f^l(FN_{\tau(i)}) -  f^l(FN_{\tau(k)}) \right).
\]
\label{thm:period_Q}
\end{thm}

The condition on the length $\ell(\alpha_i)$ of $\alpha_i$ in \textbf{Theorem \ref{thm:period_Q}} can always be fulfilled by \cite{sch}, \textbf{Corollary 4.1}. This choice is made for two reasons. First it facilitates the calculation of the length of a suitable $\alpha_{\tau(i)}$ and $\alpha_{i\tau(i)}$. Second it follows from the collar lemma in hyperbolic geometry that small simple closed geodesics have large collars, which in return gives good estimates for the upper bounds on the energies.\\
In Section 4.1 we will show how to calculate all necessary Fenchel-Nielsen coordinates. In Sections 4.2 and 4.3, we develop the functions $f^u$ and $f^l$ explicitly. In Section 4.4, we summarize these formulas and give a summary of our estimates in \textit{Table 1}. Lastly, we give a good example for our estimates in \textbf{Example 4.3}.

\subsection{Conversion of Fenchel-Nielsen coordinates for a Q-piece}

\begin{lem} Let $\mathcal{Q}_1$ be a Q-piece given in the Fenchel-Nielsen coordinates $(\ell(\beta_1),\ell(\alpha_1),{\bf tw}_1)$, where
\begin{itemize}
\item[-] $\beta_1$ is the boundary geodesic
\item[-] $\alpha_1$ an interior simple closed geodesic, such that
$\cosh(\frac{\ell(\alpha_1)}{2}) \leq \cosh(\frac{\ell(\beta_1)}{6}) + \frac{1}{2} $
\item[-] ${\bf tw}_1$ the twist parameter at $\alpha_1$.
\end{itemize}
Then there is a simple closed geodesic $\alpha_2 \subset \mathcal{Q}_1$ and a simple closed geodesic $\alpha_{12}$ in
the free homotopy class of $\alpha_1(\alpha_{2})^{-1}$, such that
\[
   \cosh(\frac{\ell(\alpha_k)}{2})= \cosh\left( \frac{\ell(\alpha_1)|t_k|}{2}\right)\sqrt{\left(\frac{\cosh( \frac{\ell(\beta_1)}{4}) }{\sinh(\frac{\ell(\ao)}{2})}\right)^2+1} , \text{ \ where \ } |t_k| =
\left\{ {\begin{array}{*{20}c}
   {|{\bf tw}_1|}\\
   {1-|{\bf tw}_1|}\\
\end{array}} \right. \text{ if } \begin{array}{*{20}c}
   {k=2}\\
   {k=12}.\\
\end{array}
\]
Furthermore, for $k \in \{2,12\}$, let ${\bf tw}_k$ be the twist parameter at $\alpha_k$, then
\[
  |{\bf tw}_k| = \min\{ \frac{2 r_k}{\ell(\alpha_k)} , 1-\frac{2 r_k}{\ell(\alpha_k)} \}, \text{ \ where \ }  r_k= \arctanh\left(\frac{\tanh(\frac{\ell(\alpha_1) |{\bf tw}_1|}{2}) \tanh(\frac{\ell(\alpha_1)}{2})}{\tanh(\frac{\ell(\alpha_k)}{2})}\right).
\]
\label{thm:FNj}
\end{lem}

\begin{figure}[h!]
\SetLabels
\L(.29*.99) $\ao'$\\
\L(.69*.99) $\ao^*$\\
\L(.34*.72) $s'$\\
\L(.38*.61) $\theta$\\
\L(.62*.35) $\theta$\\
\L(.34*.56) $r_1$\\
\L(.59*.29) $r_2$\\
\L(.45*.42) $\eta_{1}'$\\
\L(.59*.60) $\eta_{2}'$\\
\L(.52*.36) $\at'$\\
\L(.66*.25) $s^*$\\
\endSetLabels
\AffixLabels{%
\centerline{%
\includegraphics[height=5cm,width=6cm]{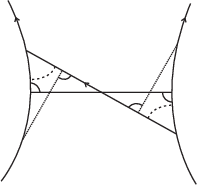}}}
\caption{Two lifts of $\ao$ in the universal covering}
\label{fig:collar_a1}
\end{figure}
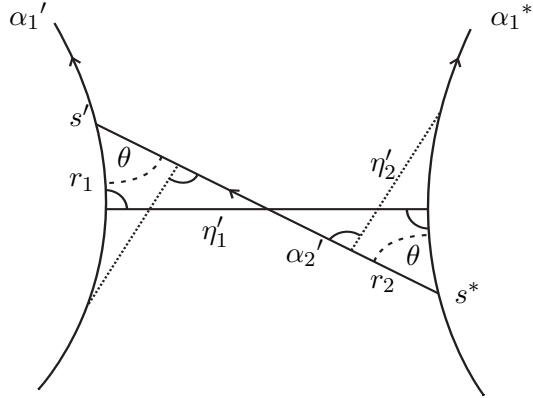
\begin{proof} In $\mathcal{Q}_1$ there exists a unique shortest geodesic arc $\eta_1$ meeting $\alpha_1$ perpendicularly on both sides of $\ao$. \textit{Figure \ref{fig:collar_a1}} shows a lift of $\ao$ and $\eta_1$ in the universal covering, $\ao$ lifts to $\ao'$ and $\ao^*$ and $\eta_1$ to $\eta_1'$. Note that $\ao'$ and $\ao^*$ have the same orientations with respect to $\eta_1'$. In the covering there exist two points, $s' \in \ao'$ and $s^* \in \ao^*$, on opposite sides of $\eta_1'$ and at the same distance from $\eta_1'$, such that $s'$ and $s^*$ are mapped to the same point $s \in \alpha_1$ by the covering map. Observe that $s' $ and $s^*$ can always be found, such that the distance $r_1$ from $\eta_1'$ is equal to $\frac{\ell(\ao) \cdot |{\bf tw}_1|}{2}$. Let $\at'$ denote the geodesic from $s'$ to $s^*$. Using $\at'$ we obtain two isometric right-angled geodesic triangles. Since $\at'$ intersects $s'$ and $s^*$ under the same angle $\theta$, the image $\at$ of $\at'$ under the universal covering map is a smooth simple closed geodesic, which intersects $\ao$ exactly once. Hence we can incorporate $\at$ into our canonical basis for $S$. Applying the cosine formula to one of the isometric triangles (see \cite{bu} p. 454), we obtain:
\[
\cosh(\frac{\ell(\at)}{2}) = \cosh(r_1) \cosh(\frac{\ell(\eta_1)}{2}), \text{ \ \ where \ \ } r_1 = \frac{\ell(\ao) \cdot |{\bf tw}_1|}{2}.
\]
The length $\ell(\eta_1)$ of $\eta_1$ can be calculated from a decomposition of $\mathcal{Q}_1$ into a Y-piece (see Equation (\ref{eq:P11})), leading to
\[
  \sinh(\frac{\ell(\eta_1)}{2})=\frac{\cosh( \frac{\ell(\beta_1)}{4}) }{\sinh(\frac{\ell(\ao)}{2})} \text{ \ thus \ } \cosh(\frac{\ell(\eta_1)}{2})=\sqrt{\left(\frac{\cosh( \frac{\ell(\beta_1)}{4}) }{\sinh(\frac{\ell(\ao)}{2})}\right)^2+1} .
\]
For further calculations we also need the angle $\theta$. From hyperbolic geometry we obtain:
\begin{equation}
\cos(\theta) = \tanh(\frac{\ell(\ao) \cdot |{\bf tw}_1|}{2}) \coth(\frac{\ell(\at)}{2})
\label{eq:cos_rho}
\end{equation}
In $\mathcal{Q}_1$, there exists likewise a unique shortest geodesic arc $\eta_2$ meeting $\alpha_2$ perpendicularly on both sides of $\at$. This arc can be seen in \textit{Figure \ref{fig:collar_a1}}. Now $\alpha_2$ and $\alpha_1$ intersect exactly once under the angle $\theta$. Consider a right-angled triangle with sides of length $\ell(\frac{\ao}{2}), r_2$ and $\ell(\frac{\eta_2}{2})$. Here $r_2$ contains information about the twist parameter ${\bf tw}_2$ with respect to $\alpha_2$. We get:
\[
   \cos(\theta) = \tanh(r_2) \coth(\frac{\ell(\ao)}{2}).
\]
Together with Equation (\ref{eq:cos_rho}), we obtain:
\[
   \tanh(r_2)= \frac{ \tanh(\frac{\ell(\ao) \cdot |{\bf tw}_1|}{2}) \coth(\frac{\ell(\at)}{2}) }{ \coth(\frac{\ell(\ao)}{2})} \text{ \ and \ } |{\bf tw}_2| = \min ( \frac{2 r_2}{\ell(\at)}, 1-\frac{2r_2}{\ell(\at)}).
\]
We will now look for a suitable $\alpha_{12}$. Consider again the lifts of $\ao$ in \textit{Figure \ref{fig:collar_a1}}. Consider the two points, $q' \in \ao'$ and $q^* \in \ao^*$ on the opposite side of $s'$ and $s^*$ with respect to the intersection point with $\eta_1$ and at distance $\ell(\ao) - r_1$ from $\eta_1'$. $q'$ and $q^*$ are mapped to the same point $q \in \mathcal{Q}_1$ by the covering map. Connecting these points we obtain a geodesic arc $\alpha_{12}'$, which maps to a simple closed geodesic $\alpha_{12}$ in $\mathcal{Q}_1$. It follows from its intersection properties with $\ao$ and $\at$ that $\alpha_{12}$ is in the free homotopy class $\ao (\at)^{-1}$.
Using the same reasoning as for $\at$, we can find its length and the twist parameter ${\bf tw}_{12}$, which leads to \textbf{Lemma \ref{thm:FNj}}.
\end{proof}

\subsection{Upper bounds for the energy of dual harmonic forms based on Q-pieces}

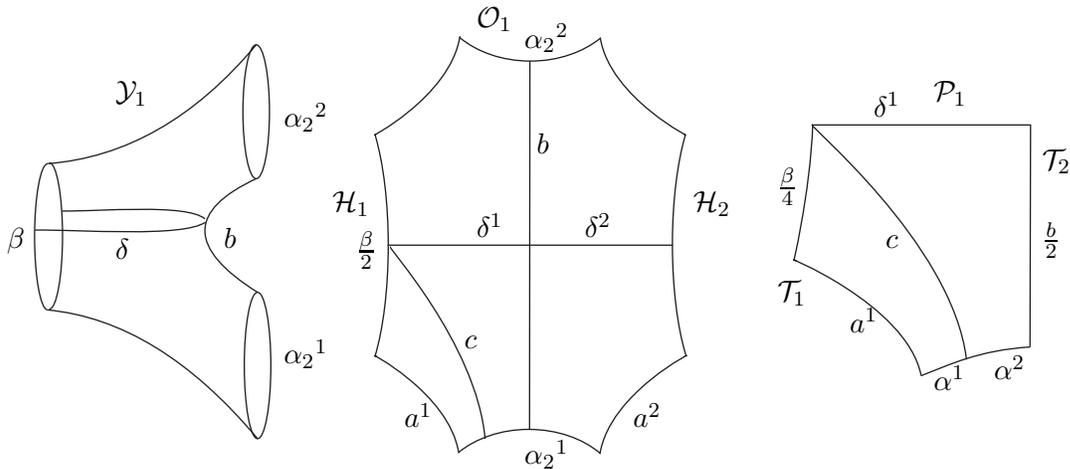
\begin{figure}[h!]
\SetLabels
\L(.13*.78) $\mathcal{Y}_1$\\
\L(.04*.50) $\beta$\\
\L(.13*.48) $\delta$\\
\L(.27*.73) $\at^2$\\
\L(.27*.28) $\at^1$\\
\L(.22*.50) $b$\\
\L(.43*.92) $\mathcal{O}_1$\\
\L(.31*.57) $\mathcal{H}_1$\\
\L(.61*.57) $\mathcal{H}_2$\\
\L(.33*.48) $\frac{\beta}{2}$\\
\L(.43*.52) $\delta^1$\\
\L(.52*.52) $\delta^2$\\
\L(.47*.88) $\at^2$\\
\L(.47*.09) $\at^1$\\
\L(.48*.68) $b$\\
\L(.37*.16) $a^1$\\
\L(.56*.16) $a^2$\\
\L(.42*.31) $c$\\
\L(.81*.78) $\mathcal{P}_1$\\
\L(.68*.60) $\frac{\beta}{4}$\\
\L(.90*.50) $\frac{b}{2}$\\
\L(.74*.34) $a^1$\\
\L(.81*.23) $\alpha^1$\\
\L(.86*.25) $\alpha^2$\\
\L(.77*.50) $c$\\
\L(.76*.75) $\delta^1$\\
\L(.68*.40) $\mathcal{T}_1$\\
\L(.90*.65) $\mathcal{T}_2$\\
\endSetLabels
\AffixLabels{%
\centerline{%
\includegraphics[height=7cm,width=15cm]{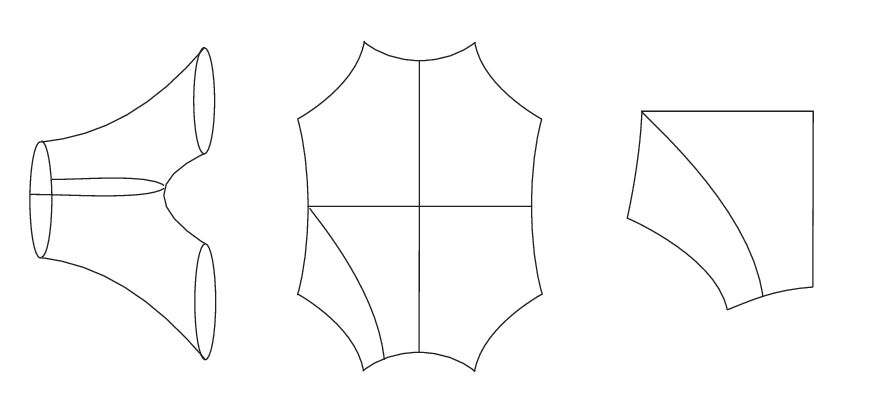}}}
\caption{Decomposition of $\mathcal{Y}_1$ into isometric hexagons $\mathcal{H}_1$ and $\mathcal{H}_2$}
\label{fig:H1H2decomp}
\end{figure}

We will establish estimates for all entries of the period Gram matrix based on the geometry of the Q-pieces $(\mathcal{Q}_i)_{i=1,3,...,2g-1}$. Following the approach given in Section 3, it is sufficient to construct suitable functions on
\[
S_{\gamma} \cap \mathcal{Q}_i,   \text{ \ where \ } \gamma \in \{\alpha_i,\alpha_{\tau(i)},\alpha_{i\tau(i)}\}, \text{ \ for \ } i \in \{1,3,...,2g-1\}.
\]
In this and the following subsection, we will only show how to obtain estimates for $E(\sigma_1)=p_{11}$ based on the geometry of $\mathcal{Q}_1$. These estimates will only depend on the Fenchel-Nielsen coordinates $(\ell(\beta_1),\ell(\at),{\bf tw}_2)$ of $\mathcal{Q}_1$. In the same way, we obtain estimates for $E(\sigma_2)=p_{22}$  based on the coordinates $(\ell(\beta_1),\ell(\ao),{\bf tw}_1)$, and for $E(\sigma_1 + \sigma_2)$ based on the coordinates $(\ell(\beta_1),\ell(\alpha_{12}),{\bf tw}_{12})$.\\
Proceeding the same way on the remaining Q-pieces and combining these estimates as described in Section 3.2 (see \textbf{Theorem \ref{thm:diag}} and \textbf{\ref{thm:ndiag}}) we finally obtain estimates for all entries of the period matrix.\\
To obtain an upper bound for $p_{11}$, we embed $S_{\at} \cap \mathcal{Q}_1$ into a hyperbolic cylinder $C$ with baseline $\at$ and denote this embedding by the same name. To obtain an estimate on $E(\sigma_1)$, we will give a parametrization of
\[
S_{\at} \cap \mathcal{Q}_1 \subset C
\]
based on a decomposition into trirectangles. To obtain this parametrization, we first cut open $\mathcal{Q}_1$ along $\at$ to obtain the Y-piece $\mathcal{Y}_1$ with boundary geodesics $\beta=\beta_1$, $\at^1$ and $\at^2$. Both $\at^1$ and $\at^2$ have length $\ell(\at)$ (see \textit{Figure \ref{fig:H1H2decomp}}).\\
Denote by $b$ the shortest geodesic arc connecting $\at^1$ and $\at^2$. We cut open $\mathcal{Y}_1$ along the shortest geodesic arcs connecting $\beta$ and the other two boundary geodesics. We call $\mathcal{O}_1$ the octagon, which we obtain by cutting open $\mathcal{Y}_1$ along these lines. By abuse of notation, we denote the geodesic arcs in $\mathcal{O}_1$ by the same letter as in $\mathcal{Y}_1$. The geodesic arc $b$ divides $\mathcal{O}_1$ into two isometric hexagons $\mathcal{H}_1$ and $\mathcal{H}_2$. This decomposition is also shown in \textit{Figure \ref{fig:H1H2decomp}}.

In $\mathcal{H}_1$ $b$ is the boundary geodesic connecting $\frac{\at^1}{2}$ and $\frac{\at^2}{2}$. Denote by $\delta^1$ the shortest geodesic arc in $\mathcal{H}_1$ connecting $b$ and the side opposite of $b$ of length $\ell(\frac{\beta}{2})$. By abuse of notation, we denote this side by $\frac{\beta}{2}$. We denote by $\delta^2$ the arc in $\mathcal{H}_2$ corresponding to $\delta^1$ in $\mathcal{H}_1$. Let $\delta=\delta^1\cup \delta^2$ be the geodesic arc in $\mathcal{O}_1$ formed by $\delta^1$ and $\delta^2$. By abuse of notation, we denote the corresponding arc in  $\mathcal{Q}_1$ and $\mathcal{Y}_1$ that maps to $\delta^1 \cup \delta^2$ in $\mathcal{O}_1$ also by $\delta$. It follows from the symmetry of $\mathcal{Y}_1$ that $\delta$ constitutes the intersection of the cut locus of $\at$ with $\mathcal{Q}_1$:
\[
   \delta = CL(\at) \cap \mathcal{Q}_1.
\]
Let $a^1$ denote the geodesic arc connecting $\frac{\at^1}{2}$ and $\frac{\beta}{2}$ in $\mathcal{H}_1$, and $a^2$ the corresponding arc in $\mathcal{H}_2$ of the same length $\ell(a^1)=\ell(a^2)=a$. Then $\delta^1$ divides $\mathcal{H}_1$ into two isometric right-angled pentagons $\mathcal{P}_1$ and $\mathcal{P}_2$. Let $\mathcal{P}_1$ be the pentagon that has $\frac{\at^1}{2}$ as a boundary. To establish the parametrization for $S_{\at} \cap \mathcal{Q}_1$, we divide $\mathcal{P}_1$ into two trirectangles. Let $c$ be the geodesic arc in $\mathcal{P}_1$ that emanates from the vertex, where $\frac{\beta}{2}$ and $\delta^1$ intersect and that meets $\frac{\at^1}{2}$ perpendicularly. It divides $\frac{\at^1}{2}$ into two parts, $\alpha^1$ and $\alpha^2$ (see \textit{Figure \ref{fig:H1H2decomp}}). $c$ divides $\mathcal{P}_1$ into two trirectangles $\mathcal{T}_1$ and $\mathcal{T}_2$, that have boundaries $\alpha^1$ and $\alpha^2$, respectively.\\
To obtain an upper bound for $p_{11}$, we need to know the geometry of $\mathcal{T}_1$ and $\mathcal{T}_2$. Hence, we need to know the lengths $a$, $\ell(\alpha^1)$, $\ell(\alpha^2)$, and $\frac{\ell(b)}{2}$. In the following subsection we will also need the length $\ell(c)$ of $c$, which we will calculate here. To obtain these lengths, we will use the geometry of $\mathcal{H}_1$, $\mathcal{P}_1$, $\mathcal{T}_1$ and $\mathcal{T}_2$. All formulas for the geometry of these polygons can be found in \cite{bu}, p. 454. From the geometry of the hyperbolic pentagon $\mathcal{P}_1$ we have:
\begin{eqnarray}
\label{eq:P11}
\sinh(\frac{\ell(b)}{2})& = & \frac{\cosh(\frac{\ell(\beta)}{4})}{\sinh(\frac{\ell(\at)}{2})} \\
\label{eq:P111}
\cosh(\ell(\delta^1))& = & \sinh(\frac{\ell(\at)}{2})\sinh(a).
\end{eqnarray}

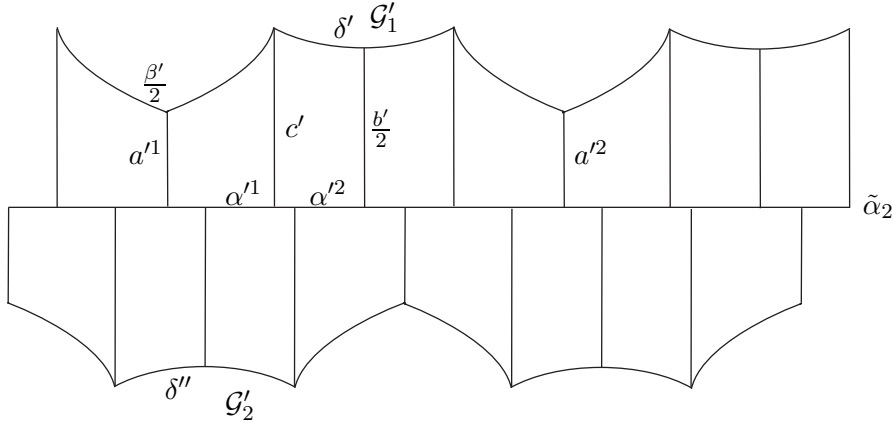
\begin{figure}[h!]
\SetLabels
\L(.45*.91) $\mathcal{G}_1'$\\
\L(.33*.05) $\mathcal{G}_2'$\\
\L(.42*.88) $\delta'$\\
\L(.28*.08) $\delta''$\\
\L(.33*.51) $\alpha'^1$\\
\L(.40*.51) $\alpha'^2$\\
\L(.26*.76) $\frac{\beta'}{2}$\\
\L(.45*.66) $\frac{b'}{2}$\\
\L(.38*.66) $c'$\\
\L(.25*.60) $a'^1$\\
\L(.62*.60) $a'^2$\\
\L(.86*.49) $\tilde{\alpha}_2$\\
\endSetLabels
\AffixLabels{%
\centerline{%
\includegraphics[height=6cm,width=12cm]{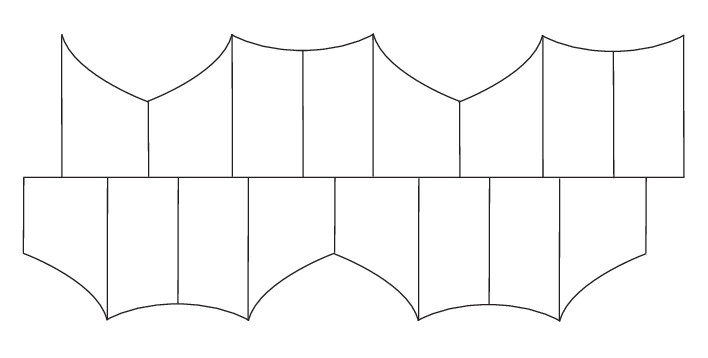}}}
\caption{Lift of $\mathcal{C}_2$ into the universal covering}
\label{fig:Y_twist}
\end{figure}

Hence, we can express $\ell(b)$ in terms of $\ell(\at)$ and $\ell(\beta)$. We obtain $a$, in terms of $\ell(b)$ and $\ell(\alpha_2)$, from the geometry of the hyperbolic hexagon $\mathcal{H}_1$ and $\ell(\delta^1) = \frac{\ell(\delta)}{2}$ in terms of $a$ and $\ell(\at)$ from Equation (\ref{eq:P111}). Finally, we can express $\ell(\alpha^2)$ and $\ell(c)$ in terms of $\ell(\delta^1)$ and $\ell(\frac{\ell(b)}{2})$ using the geometry of the hyperbolic trirectangles $\mathcal{T}_1$ and $\mathcal{T}_2$. In total, we can express the lengths $\ell(b), a, \ell(\alpha^2)$ and $\ell(\alpha^1) = \frac{\ell(\at)}{2}-\ell(\alpha^2) $ in terms of $\ell(\at)$ and $\ell(\beta)$. These formulas are simplified and summarized in Equations (\ref{eq:gon_sum1})-(\ref{eq:gon_sum3}).\\
With these formulas we can obtain a description of the boundary of $S_{\at} \cap \mathcal{Q}_1 \subset C$. Consider now $\delta \subset \mathcal{O}_1$. $\delta$ divides $\mathcal{O}_1$ into two isometric hexagons. Let $\mathcal{G}_1$ be the hexagon that contains $\at^1$ as a boundary geodesic and $\mathcal{G}_2$ be the hexagon that contains $\at^2$ as a boundary geodesic. $\delta$ forms the cut locus of $\at$ in $\mathcal{Q}_1$. Denote by $\mathcal{C}_2$ the surface that we obtain if we cut open $\mathcal{Q}_1$ along $\delta$. $\mathcal{C}_2$ is a topological cylinder around $\at$. A lift of $\mathcal{C}_2$ in the universal covering is depicted in \textit{Figure \ref{fig:Y_twist}}.

Let $\mathcal{G}_1'$ and $\mathcal{G}_2'$ denote two hexagons in this lift, that are isometric to the hexagons $\mathcal{G}_1$ and $\mathcal{G}_2$ in $\mathcal{O}_1$, and that are adjacent along the lift $\tilde{\alpha}_2$ of $\at$.  We denote by $\delta' \subset \mathcal{G}_1'$ and $\delta'' \subset \mathcal{G}_2'$ the two sides corresponding to $\delta$ in $\mathcal{O}_1$. We keep the notation from $\mathcal{O}_1$, but denote all corresponding geodesic arcs in the covering space with prime, i.e. a lift of $\alpha^1$ is denoted by $\alpha'^1$ etc.\\
In the lift of $\mathcal{C}_2$ the two hexagons $\mathcal{G}_1'$ and $\mathcal{G}_2'$ are shifted against each other by the length $|{\bf tw}_2| \cdot \ell(\at)$. It can be seen from \textit{Figure \ref{fig:Y_twist}}, how to parametrize $S_{\at} \cap \mathcal{Q}_1$ in a cylinder $C$ around $\at$. Here all boundaries are boundaries of trirectangles, which are isometric to either $\mathcal{T}_1$ or $\mathcal{T}_2$, which can be parametrized in Fermi coordinates.  Using these formulas in \textbf{Theorem \ref{thm:capa_S2}}, we can find an upper bound $f^u(FN_2)$ for $E(\sigma_1)$:
\[
  f^u(FN_2) \geq \capa(S_{\at} \cap \mathcal{Q}_1) \geq E(\sigma_1) =p_{11}.
\]
We obtain a simplified upper bound, if we define our test function only on the collar $Z_{\min\{a,\frac{\ell(b)}{2}\}}(\at)$ (see defintion (\ref{eq:zrx})). This upper bound $f^u_{simp}$ corresponds to the method from \cite{bs} applied to a Q-piece and is given in inequality (\ref{flu_simp}). These formulas are summarized in Section 4.4 and the results are summarized in \textit{Table 1}.

\subsection{Lower bounds for the energy of dual harmonic forms based on Q-pieces}\label{sec:lowQ}

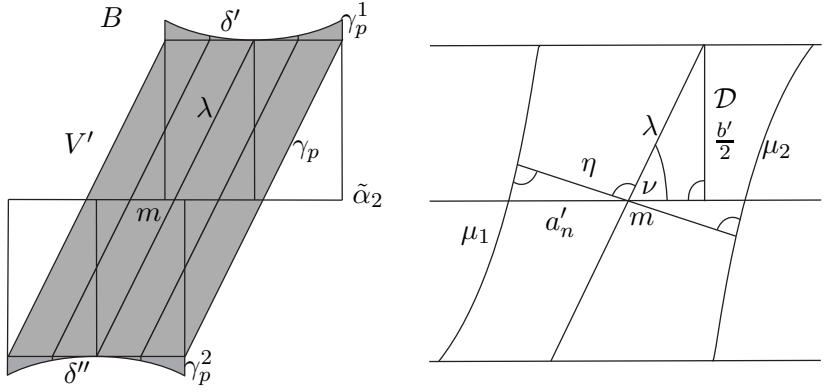
\begin{figure}[h!]
\SetLabels
\L(.23*.88) $B$\\
\L(.20*.60) $V'$\\
\L(.31*.68) $\lambda$\\
\L(.26*.45) $m$\\
\L(.33*.87) $\delta'$\\
\L(.20*.09) $\delta''$\\
\L(.44*.49) $\tilde{\alpha}_2$\\
\L(.30*.11) $\gamma_p^2$\\
\L(.39*.60) $\gamma_p$\\
\L(.43*.88) $\gamma_p^1$\\
\L(.74*.70) $\mathcal{D}$\\
\L(.53*.41) $\mu_1$\\
\L(.78*.60) $\mu_2$\\
\L(.63*.56) $\eta$\\
\L(.68*.51) $\nu$\\
\L(.67*.44) $m$\\
\L(.68*.64) $\lambda$\\
\L(.60*.43) $a_n'$\\
\L(.74*.61) $\frac{b'}{2}$\\
\endSetLabels
\AffixLabels{%
\centerline{%
\includegraphics[height=6cm,width=12cm]{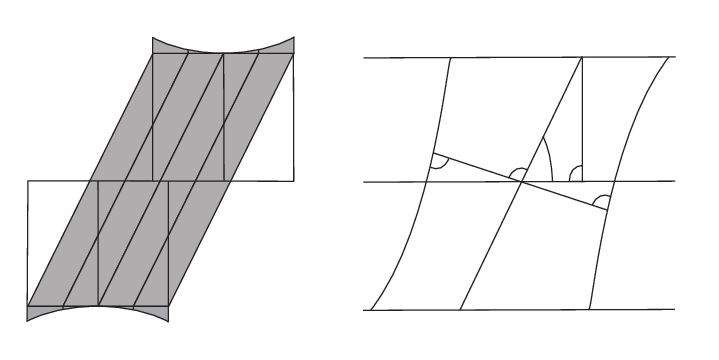}}}
\caption{The area $B$ \textit{(grey)} and the construction of skewed Fermi coordinates $\psi^{\nu}$}
\label{fig:Fermi_nu}
\end{figure}

Consider a primitive $F_1$ of $\sigma_1$ in $\mathcal{C}_2=S_{\at} \cap \mathcal{Q}_1 \subset C$. The two geodesic arcs $\delta'$ and $\delta''$ corresponding to $\delta \subset \mathcal{Q}_1$ constitute $CL(\at)^{red} \cap \partial \mathcal{C}_2$. We will use the theoretical approach from Section 3 to obtain a concrete lower bound $f^l(FN_2)$ for
\[
      p_{11} = E_S(F_1) > E_{B}(F_1) \geq f^l(FN_2) , \text{ \ where \ } B=S_2^{red} \cap  \mathcal{C}_2.
\]
We will give a suitable construction for $B=S_2^{red} \cap  \mathcal{C}_2$ in Section 4.3.1. To this end, we lift  $\mathcal{C}_2$ into the universal covering as in the previous subsection (see \textit{Figure \ref{fig:Y_twist}}). We use the same notation for the geodesic arcs that occur. The important cut-out from \textit{Figure \ref{fig:Y_twist}} is depicted in \textit{Figure \ref{fig:Fermi_nu}}.\\

Let $B$ be the grey hatched subset in the lift of $\mathcal{C}_2$ in \textit{Figure \ref{fig:Fermi_nu}}. We will now give an exact description and parametrization of $B$.

\subsubsection{Parametrization of $B=S_2^{red} \cap  \mathcal{C}_2$}

The boundary of $B$ contains the lines $\delta'$ and $\delta''$. For each point $p_1 \in \delta'$, there exists a point $p_2 \in \delta''$, such that $p_1$ and $p_2$ map to the same point $p$ on $\delta \subset \mathcal{Q}_1$. We may assume, without loss of generality, that
\[
    F_1(p_2)-F_1(p_1)=1 \text { \ for all \ } p_1 \in \delta'.
\]
We will describe $B$ as a union of lines, where each line $l_{p}$ connects $p_1$ and $p_2$. The line $l_{p}$ is defined as follows. From $p_1$ we go along the geodesic that meets $\tilde{\alpha}_2$ perpendicularly until we meet $\partial Z_{\frac{\ell(b)}{2}}(\tilde{\alpha}_2)$, the boundary of the collar $Z_{\frac{\ell(b)}{2}}(\tilde{\alpha}_2)$ (see definition (\ref{eq:zrx})). We call this intersection point $p_1'$ and the geodesic arc that forms $\gamma_p^1$. Let $p_2'$ be the point on $\partial Z_{\frac{\ell(b)}{2}}(\tilde{\alpha}_2)$ on the other side of $\tilde{\alpha}_2$ that can be reached analogously, starting from $p_2$. We now go along the geodesic arc that connects $p_1'$ and $p_2'$. We call this arc $\gamma_{p}$. Then from $p_2'$, we move along the geodesic arc connecting $p_2'$ and $p_2$. We call this arc $\gamma_p^2$. We define $l_{p}$ as the line traversed in this way. Let $B$ be the disjoint union of these lines:
\[
     B= \biguplus\limits_{p \in \delta} \{l_{p}\}
\]
Let $\lambda$ be the geodesic arc connecting the midpoints of $\delta'$ and $\delta''$, and let $m$ be the midpoint  of $\lambda$.
We will use a bijective parametrization  $\varphi:(t,s) \mapsto \varphi(t,s)$ of $B$, such that
\begin{itemize}
\item[-] $ \varphi(0,0) = m$
\item[-] for all $t \in [-\ell(\alpha^2),\ell(\alpha^2)]$, $ \varphi(t,0) \in \tilde{\alpha}_2$ has directed distance $t$ from $m$
\item[-] for a fixed $t_0 \in [-\ell(\alpha^2),\ell(\alpha^2)]$, $\varphi(t_0,\cdot)$ parametrizes the line $l_{p}$ that traverses $\tilde{\alpha}_2$ in a point with directed distance $t_0$ from $m$ by arc length.
\end{itemize}
We parametrize the sets $\bigcup\limits_{p \in \delta} \{\gamma_p^1\}$ and $\bigcup\limits_{p \in \delta} \{\gamma_p^2\}$ in Fermi coordinates with baseline $\tilde{\alpha}_2$. The proper parametrization can be deduced from the geometry of the trirectangle $\mathcal{T}_2$.\\
\\
We will parametrize $Z_{\frac{\ell(b)}{2}}(\tilde{\alpha}_2) \cap B  = \bigcup\limits_{p \in \delta} \{\gamma_p\}$ using \textit{skewed Fermi coordinates} $\psi^{\nu}$, with angle $\nu$ and baseline $\tilde{\alpha}_2$. These are defined in the same way as the usual Fermi coordinates $\psi$ (see Section 2.1), but instead of moving along geodesics emanating perpendicularly from the baseline, we move along geodesics that meet the baseline under the angle $\nu$. We will not give these coordinates explicitly, but will derive the essential information from the Fermi coordinates $\psi$.\\
We remind the reader that $\lambda$ is the geodesic arc connecting the midpoints of $\delta'$ and $\delta''$. Its midpoint $m$ and the endpoints of $\frac{b'}{2}$ are the vertices of a right-angled triangle $\mathcal{D}$ (see \textit{Figure \ref{fig:Fermi_nu}}). In our case the angle $\nu$ for the coordinates $\psi^{\nu}$ is the angle of $\mathcal{D}$ at the midpoint $m$. It follows from the geometry of right-angled triangles that
\begin{equation}
\cosh(\frac{\ell(\lambda)}{2})=\cosh(\frac{\ell(b)}{2})\cosh(\frac{\ell(\at) |{\bf tw}_2|}{2}),
\label{eq:twist_D}
\end{equation}
where we assume, without loss of generality, that the twist parameter ${\bf tw}_2$ is in the interval $[0,\frac{1}{2}]$. Otherwise the situation is symmetric to the depicted one. Using the geometry of the right-angled triangle $\mathcal{D}$ we have:
\begin{equation}
\sin(\nu) = \frac{\sinh(\frac{\ell(b)}{2})}{\sqrt{\cosh(\frac{\ell(b)}{2})^2 \cosh(\frac{\ell(\at) |{\bf tw}_2|}{2})^2-1}}.
\label{eq:twist_nu}
\end{equation}
Consider the following geodesic arcs in $Z_{\frac{\ell(b)}{2}}(\tilde{\alpha}_2) \cap B$. For a $n \in \N$, let $a_n'$ be a geodesic arc of length $\frac{2\alpha^2}{n}$ on $\tilde{\alpha}_2$ with midpoint $m$. $\lambda$ intersects $a_n'$ in $m$ under the angle $\nu$. This is depicted in \textit{Figure \ref{fig:Fermi_nu}}.

Let $\eta'$ be a geodesic intersecting $\lambda$ perpendicularly in $m$. Let $\mu_1$ and $\mu_2$ be two geodesic arcs with endpoints on $Z_{\frac{\ell(b)}{2}}(\tilde{\alpha}_2)$ that intersect $\eta'$ perpendicularly, such that each of the arcs passes through an endpoint of $a_n'$ on each side of $\lambda$. Let $\eta$ be the geodesic arc on $\eta'$ with endpoints on $\mu_1$ and $\mu_2$. For fixed $n \in \N$, we denote by $\eta_n$ the length of $\eta$ and by $\mu^n$ the length of $\mu_1$ and $\mu_2$:
\[
 \eta_n = \ell(\eta) \text{ \ \ and \ \ }  \mu^n = \ell(\mu_1)=\ell(\mu_2).
\]
By choosing usual Fermi coordinates with baseline $\eta$, we can parametrize the strip, whose boundary lines are $\mu_1$ and $\mu_2$ and two segments of $\partial Z_{\frac{\ell(b)}{2}}(\tilde{\alpha}_2)$ (see \textit{Figure \ref{fig:Fermi_nu}}).\\
$n$ such strips can be aligned next to each other to obtain a parametrization of $Z_{\frac{\ell(b)}{2}}(\tilde{\alpha}_2) \cap B$. For $n \rightarrow \infty$ we obtain a parametrization $\psi^{\nu}$ of $Z_{\frac{\ell(b)}{2}}(\tilde{\alpha}_2) \cap B$. We get:
\begin{equation}
     \lim\limits_{n \to \infty} n \cdot \eta_n = \sin(\nu)2\ell(\alpha^2) \text{ \ \ and \ \ }  \lim\limits_{n \to \infty} \mu^n = \ell(\lambda).
\label{eq:ell_eta} 
\end{equation}
Combining the parametrizations for the several pieces of $B$, we may assume that we have a parametrization $\varphi$ that satisfies our conditions.

\subsubsection{Evaluating the lower bound for  $p_{11} = E_S(F_1)$} \label{sec:low_expl}

Let $\varphi(\{t_0\} \times [-x,x])$ be the parametrization for a line $l_p = \gamma_p^1 \cup \gamma_p \cup \gamma_p^2 \subset B$, such that  
\[
p_1=\varphi(t_0,-x)  \in \delta' \text{ \  and \ } p_2=\varphi(t_0,x) \in \delta''. 
\]
The function $F_1$ satisfies the boundary conditions $F_1(p_2)=1+\tilde{c}$ and $F_1(p_1)=\tilde{c}$, where $\tilde{c}$ is a constant. As our estimate depends only on the difference $F_1(p_2)-F_1(p_1) =1$, the constant $\tilde{c}$ is not important for our estimate and we assume that $\tilde{c}=0$. \\
The lower bound for  $E_S(F_1)$ is obtained by projecting the tangent vectors of $F_1$ onto the curves  $(\varphi(\{t_0\} \times [-x,x]))_{t_0}$ of our parametrization. This can be seen as a limit process, where we consider this projection on aligned  strips $(\varphi([t_0-\epsilon,t_0]\times[-x,x])$ as described in the previous subsection. The limit, however, does not depend on the width of the strip and it is therefore sufficient to consider single lines, which we will do in the following.\\  
We know that by the definition of $\frac{\lambda}{2}$:
\begin{eqnarray*}
F_1(\varphi(t_0,-x))&=&0 \text{ \ \ and \ \ } F_1(\varphi(t_0,x))=1 \text{ \ \ and \ \ } \\
F_1(\varphi(t_0,\frac{-\ell(\lambda)}{2}))&=&a_1 \text{ \ and \ \ } F_1(\varphi(t_0,\frac{\ell(\lambda)}{2}))=a_2.\\
\end{eqnarray*}
We first focus on the second condition above for the boundary of the geodesic segment $\varphi(t_0\times[\frac{-\ell(\lambda)}{2},\frac{\ell(\lambda)}{2}])$. We can consider skewed Fermi coordinates $\psi^{\nu}$ as a limit case of Fermi coordinates with respect to an imaginary baseline $\eta$ (see \textit{Figure \ref{fig:Fermi_nu}}). Let $F_{t_0}=f_{t_0} \circ (\psi^{\nu})^{-1}$ be a function defined on $\varphi(\{t_0\} \times[\frac{-\ell(\lambda)}{2},\frac{\ell(\lambda)}{2}])$, such that $f_{t_0}$ realizes the minimum
\[
\min \{\sin(\nu) \cdot \int_{\frac{-\ell(\lambda)}{2}}^{\frac{\ell(\lambda)}{2}} \cosh(s)f'(s)^2 \, ds \mid f \in \lip([\frac{-\ell(\lambda)}{2},\frac{\ell(\lambda)}{2}])), f(\frac{-\ell(\lambda)}{2})=a_1 \text{ and } f(\frac{\ell(\lambda)}{2})=a_2\}  \text{ \ (see (\ref{eq:projection})).}
\]
Comparing with Equation (\ref{eq:energy}), we see that $F_{t_0}$ is the minimizing function for the projection of tangent vectors onto $\varphi(\{t_0\} \times[\frac{-\ell(\lambda)}{2},\frac{\ell(\lambda)}{2}])$. Here the correction factor $\sin(\nu)$ follows the fact that our baseline should be orthogonal to $\lambda$ (see (\ref{eq:ell_eta})). By applying the calculus of variations (see \cite{ge}, p. 14-16) to the above integral, we obtain analogously to the construction for \textbf{Theorem \ref{thm:capa_S2}} that
\begin{equation}
\sin(\nu) \cdot \int_{\frac{-\ell(\lambda)}{2}}^{\frac{\ell(\lambda)}{2}} \cosh(s)f_{t_0}'(s)^2 \,ds = \frac{(a_2-a_1)^2 \sin(\nu)}{2(\arctan(\exp(\frac{\ell(\lambda)}{2}))-\arctan(\exp(-\frac{\ell(\lambda)}{2})))}=k_1(a_2-a_1)^2.
\label{eq:low_V1}
\end{equation}
We can extend $F_{t_0}$ to a function on $\varphi(\{t_0\} \times[-x,x])$ that satisfies the boundary conditions of $F_1$. As before, we choose $F_{t_0}$ such that it minimizes the projection of tangent vectors onto the two disjoint geodesic segments $\varphi(\{t_0\} \times[-x,\frac{-\ell(\lambda)}{2}])$ and $\varphi(\{t_0\} \times[\frac{\ell(\lambda)}{2},x])$. Using the parametrization in Fermi coordinates $F_{t_0} = f_{t_0} \circ \psi^{-1}$ on these segments, we get:
\begin{eqnarray}
\int_{[- x(t_0),-\frac{\ell(b)}{2}] \cup [\frac{\ell(b)}{2}, x(t_0)] } \cosh(s)f_{t_0}'(s)^2 \,ds= \frac{a_1^2+(1-a_2)^2 }{2(\arctan(\exp( x(t_0)))-\arctan(\exp(\frac{\ell(b)}{2})))}=  \nonumber \\ 
  k_2(t_0)(a_1^2+(1-a_2)^2), \text{ \ where \ }  x(t_0)=\ell(\gamma_p^1) + \frac{\ell(b)}{2}. 
\label{eq:low_V2}
\end{eqnarray}
As $a_1=F_1(\varphi(t_0,\frac{-\ell(\lambda)}{2}))$ and $a_2=F_1(\varphi(t_0,\frac{\ell(\lambda)}{2}))$, we have by construction $E_{l_p}(F_1) \geq  E(F_{t_0})$.\\
Though we do not know the values $a_1$ and $a_2$, we obtain a lower bound of the energy of $F_1$, if we determine the values $F_{t_0}(\varphi(t_0,\frac{-\ell(\lambda)}{2}))=c_1=c_1(t_0)$ and $F_{t_0}(\varphi(t_0,\frac{\ell(\lambda)}{2}))=c_2=c_2(t_0)$, respectively, such that these values are minimizing the total energy $E_{l_p}(F_{t_0})$. As the two arcs $\gamma_p^1$ and $\gamma_p^2$ have the same length, we have to solve the following problem: using Equations (\ref{eq:low_V1}) and (\ref{eq:low_V2}) we have to find $c_1,c_2$, such that $1-c_2 =c_1 \Leftrightarrow (c_2-c_1)=1-2c_1$, and
\[
     k_1(c_2-c_1)^2 + k_2(t_0)(c_1^2+(1-c_2)^2)   
\]
is minimal. We obtain that $c_1=\frac{k_1}{k_2(t_0)+2k_1}$. As $B = \biguplus\limits_{p \in \delta} \{l_{p}\}$ we obtain in total by integration:
\begin{equation}
   p_{11}= E(F_1) \geq  2 \int \limits_{t=0}^{\ell(\alpha^2)} \frac{k_1 k_2(t)}{k_2(t) +2 k_1} \,dt = f^l(\ell(\beta_1),\ell(\at),{\bf tw}_2).
\label{eq:q11_toolow}
\end{equation}
To obtain the lower bound $f^l$ that depends only on $\ell(\at), |{\bf tw}_2|$, and $\ell(\beta_1)$ we first have to express $\frac{\ell(b)}{2}$ and $\ell(\lambda)$  and $\nu$ in terms of these variables (see Equations (\ref{eq:P11}),(\ref{eq:twist_D}) and (\ref{eq:twist_nu})). Using the parametrization of $\mathcal{T}_2$ in Equation (\ref{eq:low_V2}), we can then express $f^l$ in terms of $\ell(\at), |{\bf tw}_2|$ and $\ell(\beta_1)$. This way we obtain explicit values in Equation (\ref{eq:q11_toolow}). These formulas are summarized in the following subsection.

\subsection{Summary}
In this section, we summarize the formulas from the previous subsections and outline our estimates in \textit{Table 1}. We also give an example for our estimates in \textbf{Example 4.3}. First, we give a description of $f^u$ and $f^l$ from \textbf{Theorem \ref{thm:period_Q}}.

\subsubsection{Upper bound $f^u$ from \textbf{Theorem \ref{thm:period_Q}}}

In the remaining part of this paper we fix the notation in the following way: for $j\in \{i,\tau(i),i\tau(i)\}$, let $\mathcal{Q}_i$ be a Q-piece given in Fenchel-Nielsen coordinates $FN_j=(\ell(\beta_j),\ell(\alpha_j),{\bf tw}_j)$, where $\beta_j=\beta_i$ is the boundary geodesic of $\mathcal{Q}_i$, and ${\bf tw}_j \in  (-\frac{1}{2},\frac{1}{2}]$ be the twist parameter at an interior simple closed geodesic $\alpha_j$. We have from Section 4.2:

\begin{eqnarray}
\sinh(\frac{\ell(b)}{2})&=& \frac{\cosh(\frac{\ell(\beta_j)}{4})}{\sinh(\frac{\ell(\alpha_j)}{2})} \label{eq:gon_sum1}\\ 
\coth(a) &=&\tanh(\frac{\ell(b)}{2})\cosh(\frac{\ell(\alpha_j)}{2}) \text{ and } \sinh(\ell(c))=\frac{\cosh(\frac{\ell(\beta)}{4})}{\sqrt{\tanh(\frac{\ell(b)}{2})^2\cosh(\frac{\ell(\alpha_j)}{2})^2-1}}. \label{eq:gon_sum2} \\
\coth(\ell(\alpha^2))&=&\cosh(\frac{\ell(b)}{2})^2\tanh(\frac{\ell(\alpha_j)}{2}) \text{ and }
\ell(\alpha^1)=\frac{\ell(\alpha_j)}{2}-\ell(\alpha^2).
\label{eq:gon_sum3}
\end{eqnarray}

Using the above, we obtain a description of the cut locus $CL(\alpha_j) \cap \mathcal{Q}_i$ in a cylinder $\mathcal{C}_j$ in Fermi coordinates. Set $\Sp^1_{\alpha_j} = \R \mod (t \mapsto t+\ell(\alpha_j))$. For $l \in \{1,2\}$ let
\[
a_l: \Sp^1_{\alpha_j}  \rightarrow \R ,  a_l: t \mapsto a_l(t)
\]
be a parametrization of the two connected components of $CL(\alpha_j) \cap \mathcal{Q}_i$ in $\mathcal{C}_j$. Then

\begin{eqnarray*}
a_2(t):&=&
\left\{ {\begin{array}{*{20}c}
   {\arctanh(\cosh(t-\ell(\alpha^2))\tanh(\frac{\ell(b)}{2}))}  \\
   {\arctanh(\cosh(t-(2\ell(\alpha^2)+\ell(\alpha^1)))\tanh(\frac{a}{2}))}  \\
\end{array}} \right.\text{ \ if \ } \begin{array}{*{20}c}
   { t \in(0, 2\ell(\alpha^2)]}  \\
   {t \in (2\ell(\alpha^2), \ell(\alpha_j)] }  \\
\end{array}\\
a_1(t):&=& - a_2(t + |{\bf tw}_j|)
\end{eqnarray*}

Applying \textbf{Theorem \ref{thm:capa_S2}} to estimate the capacity of $S_{\alpha_j} \cap \mathcal{Q}_i$  with boundary $CL(\alpha_j) \cap \mathcal{Q}_i$, we obtain:
\begin{eqnarray*}
f^u(FN_j):= \int\limits_{t=0}^{\ell(\alpha_j)} {\frac{1+ \frac{1}{3}\cdot\left(\frac{(a_1'(t))^2}{\cosh^2(a_1(t))}+\frac{a_1'(t)}{\cosh(a_1(t))}\cdot\frac{a_2'(t)}{\cosh(a_2(t))}+\frac{(a_2'(t))^2}{\cosh^2(a_2(t))}\right)}{2(\arctan(\exp(a_2(t)))-\arctan(\exp(a_1(t))))} \, dt }  \geq  \\ \capa(S_{\alpha_j} \cap \mathcal{Q}_i) \geq \text{ \ \ \ \ \ \ \ \ \ \ \ \ \ \ \ \ \ \ \ \ \ \ \ \ \  }  \\
\int\limits_{t=0}^{\ell(\alpha_j)} {\frac{1}{2(\arctan(\exp(a_2(t)))-\arctan(\exp(a_1(t))))} \, dt} := f^u_{low}(FN_j),
\end{eqnarray*}
where $f^u_{low}(FN_j)$ is a lower bound for the capacity of $S_{\alpha_j} \cap \mathcal{Q}_i$. For the simplified upper bound $f^u_{simp}(FN_j)$ that corresponds to the method in 
\cite{bs} we have:
\begin{equation}
  f^u_{simp}(FN_j) = \frac{\ell(\alpha_j)}{2(\arctan(e^{\min\{a,\frac{\ell(b)}{2}\}})-\arctan(e^{-\min\{a,\frac{\ell(b)}{2}\}}))}  \geq E(\sigma_{\tau(j)}).
\label{flu_simp}  
\end{equation}
This upper bound is in the range of $f^u$ for small $\alpha_j$, but in general much larger.

\subsubsection{Lower bound $f^l$ from \textbf{Theorem \ref{thm:period_Q}}}

Based on Section \ref{sec:lowQ}, we first give a suitable construction for $S_j^{red} \cap \mathcal{Q}_i$, where $j\in \{i,\tau(i),i\tau(i)\}$. From Equations (\ref{eq:twist_D}) and (\ref{eq:twist_nu}) we obtain (see \textit{Figure \ref{fig:Fermi_nu}}):
\begin{equation*}
\cosh(\frac{\ell(\lambda)}{2})=\cosh(\frac{\ell(b)}{2})\cosh(\frac{\ell(\alpha_j) |{\bf tw}_j|}{2}) \text{ \ and \ }
\sin(\nu) = \frac{\sinh(\frac{\ell(b)}{2})}{\sqrt{\cosh(\frac{\ell(b)}{2})^2 \cosh(\frac{\ell(\alpha_j) |{\bf tw}_j|}{2})^2-1}}.
\end{equation*}
Using the above we obtain a description of the cut locus $CL(\alpha_j)^{red} \cap \mathcal{Q}_i$ in a cylinder $\mathcal{C}_j$ in Fermi coordinates. Let
\[
a_{red}: [0,2\ell(\alpha^2)]  \rightarrow \R ,  a_{red}: t \mapsto a_{red}(t)
\]
be a parametrization of one of the two connected components of $CL(\alpha_j)^{red} \cap \mathcal{Q}_i$. Then

\[
a_{red}(t):= \arctanh(\cosh(t-\ell(\alpha^2))\tanh(\frac{\ell(b)}{2}))  \text{ \ for  \ } t \in [0,  2\ell(\alpha^2)].
\]

From Equation (\ref{eq:low_V1}) and (\ref{eq:low_V2}) (see Section \ref{sec:low_expl}) we have:
\begin{eqnarray*}
k_1 &=&  \frac{\sin(\nu)}{2(\arctan(\exp(\frac{\ell(\lambda)}{2}))-\arctan(\exp(-\frac{\ell(\lambda)}{2})))}\\
k_2(t):&=& \frac{1}{2(\arctan(\exp(a_{red}(t)))-\arctan(\exp(\frac{\ell(b)}{2})))} \text{ \ for \ } {t \in [0,2\ell(\alpha^2)]}.
\end{eqnarray*}

Finally, we obtain the lower bound $f^l(FN_j)$ on $E(\sigma_{\tau(j)})$, where $ \sigma_{\tau(i\tau(i))} = \sigma_i + \sigma_{\tau(i)}$ from Equation (\ref{eq:q11_toolow}):
\[
    E(\sigma_{\tau(j)}) \geq   2 \int \limits_{t=0}^{\ell(\alpha^2)} \frac{k_1 k_2(t)}{k_2(t) +2 k_1} \,dt := f^l(FN_j).
\]
From $(f^u(FN_j))_j$ and $(f^l(FN_j))_j$ all entries of $P_S$ can be estimated. This follows from \textbf{Theorem \ref{thm:diag}} and \textbf{\ref{thm:ndiag}}.\\
The following table provides a comparison of the estimates for the energy of a harmonic form based on the geometry of a Q-piece $\mathcal{Q}_i$, given in Fenchel-Nielsen coordinates $FN_j=(\ell(\beta_j),\ell(\alpha_j),{\bf tw}_j)$ for ${\bf tw}_j=0$ and ${\bf tw}_j=\frac{1}{4}$. \\

\begin{table}[htbp]
\begin{center}
\begin{tabular}{|c|c|c|c|c|c|c|c|}
\hline
\multicolumn{1}{|l}{} & \multicolumn{1}{l|}{}  & \multicolumn{1}{l}{} & \multicolumn{1}{l}{${\bf tw}_j=0$} & \multicolumn{1}{l|}{} & \multicolumn{1}{l}{} & \multicolumn{1}{l}{${\bf tw}_j=\frac{1}{4}$} & \multicolumn{1}{l|}{} \\ \hline
\multicolumn{1}{|l|}{$\ell(\beta_j)$} & \multicolumn{1}{l|}{$\ell(\alpha_j)$}  & \multicolumn{1}{l|}{$f^u(FN_j)$} & \multicolumn{1}{l|}{$f^u_{low}(FN_j)$} & \multicolumn{1}{l|}{$f^l(FN_j)$} & \multicolumn{1}{l|}{$f^u(FN_j)$} & \multicolumn{1}{l|}{$f^u_{low}(FN_j)$} & \multicolumn{1}{l|}{$f^l(FN_j)$}  \\ \hline
\multicolumn{ 1}{|c|}{} & 1  & 0.55 & 0.42 & 0.40 & 0.55 & 0.41 & 0.39  \\
\multicolumn{ 1}{|c|}{} & 2 & 1.41 & 1.14 & 1.11 & 1.43 & 1.12 & 1.00 \\
\multicolumn{ 1}{|c|}{1} & 5 & 8.70 & 8.17 & 8.13  & 7.73 & 7.00 & 0.90 \\
\multicolumn{ 1}{|c|}{} & 10 & 112.46 & 111.85 & 111.80 & 61.96 & 60.94 & 0.04  \\
\multicolumn{ 1}{|c|}{} & 20 & 16772.11 & 16771.50 & 16771.45 & 2750.28 & 2749.10 & 0.00011  \\ \hline
\multicolumn{ 1}{|c|}{} & 1 & 0.47 & 0.41 & 0.33 & 0.47 & 0.40 & 0.33 \\
\multicolumn{ 1}{|c|}{} & 2 & 1.23 & 1.08 & 0.95 & 1.22 & 1.07 & 0.87 \\
\multicolumn{ 1}{|c|}{2} & 5 & 7.87 & 7.49 & 7.30 & 6.92 & 6.44 & 0.91 \\
\multicolumn{ 1}{|c|}{} & 10 & 102.76 & 102.30 & 102.11 & 56.48 & 55.76 & 0.04 \\
\multicolumn{ 1}{|c|}{} & 20 & 15340.96 & 15340.49 & 15340.22 & 2515.41 & 2514.53 & 0.00012 \\ \hline
\multicolumn{ 1}{|c|}{} & 1 & 0.44 & 0.40 & 0.12 & 0.44 & 0.40 & 0.12 \\
\multicolumn{ 1}{|c|}{} & 2 & 1.10 & 1.01 & 0.36 & 1.10 & 1.01 & 0.34 \\
\multicolumn{ 1}{|c|}{5} & 5 & 5.41 & 5.17 & 3.82 & 5.12 & 4.82 & 0.92  \\
\multicolumn{ 1}{|c|}{} & 10 & 62.08 & 61.71 & 60.30 & 35.62 & 35.11 & 0.06  \\
\multicolumn{ 1}{|c|}{} & 20 & 9161.19 & 9160.80 & 9159.37 & 1504.60 & 1503.90 & 0.00018 \\ \hline
\multicolumn{ 1}{|c|}{} & 1 & 0.58 & 0.42 & 0.01 & 0.59 & 0.42 & 0.010 \\
\multicolumn{ 1}{|c|}{} & 2 & 1.41 & 1.14 & 0.03 & 1.42 & 1.12 & 0.032 \\
\multicolumn{ 1}{|c|}{10} & 5 & 6.41 & 6.13 & 0.65 & 5.78 & 5.46 & 0.41  \\
\multicolumn{ 1}{|c|}{} & 10 & 26.06 & 25.60 & 17.54 & 19.40 & 18.75 & 0.14  \\
\multicolumn{ 1}{|c|}{} & 20 & 2828.11 & 2827.62 & 2819.53 & 484.70 & 483.96 & 0.14 \\ \hline
\multicolumn{ 1}{|c|}{} & 1 & 0.69 & 0.42 & 0.000068 & 0.72 & 0.42 & 0.000068 \\
\multicolumn{ 1}{|c|}{} & 2 & 1.86 & 1.17 & 0.000212 & 1.83 & 1.15 & 0.000205 \\
\multicolumn{ 1}{|c|}{20} & 5 & 9.18 & 8.41 & 0.004493 & 8.32 & 7.21 & 0.003701 \\
\multicolumn{ 1}{|c|}{} & 10 & 82.11 & 81.71 & 0.66  & 45.52 & 44.96 & 0.18 \\
\multicolumn{ 1}{|c|}{} & 20 & 347.17 & 346.61 & 231.58 & 99.69 & 98.69 & 0.0042 \\ \hline
\end{tabular}
\end{center}
\caption{Comparison of the estimates for the energy of a harmonic form based on the geometry of a Q-piece $\mathcal{Q}_i$, given in Fenchel-Nielsen coordinates $FN_j=(\ell(\beta_j),\ell(\alpha_j),{\bf tw}_j)$, for $j \in \{i,\tau(i),i\tau(i)\}$ and ${\bf tw}_j \in \{0,\frac{1}{4}\}$.}
\label{tab:notwist}
\end{table}

\textbf{Example 4.3} Let $\mathcal{Q}_1$ and $\mathcal{Q}_3$ be two isometric Q-pieces given in Fenchel-Nielsen coordinates $FN_1$ and $FN_3$, respectively, where
\[
FN_i = (\ell(\beta_i),\ell(\alpha_i),{\bf tw}_i) = (2,1,0.1), \text{ \ for \ } i\in \{1,3\},
\]
where $\beta_i$ is the boundary geodesic, $\alpha_i$ an interior simple closed geodesic, and ${\bf tw}_i$ the twist parameter at $\alpha_i$. Let
\[
S=\mathcal{Q}_1+\mathcal{Q}_3
\]
be a Riemann surface of genus $2$, which we obtain by gluing $\mathcal{Q}_1$ and $\mathcal{Q}_3$ along $\beta_1$ and $\beta_3$ with arbitrary twist parameter ${\bf tw}_{\beta} \in (-\frac{1}{2},\frac{1}{2}]$. Then there exists a canonical basis ${\rm A}= (\alpha_1,\alpha_2,\alpha_3,\alpha_4)$ and a corresponding period Gram matrix $P_S$, such that
\[
\left( {\begin{array}{*{20}c}
 2.11 & -0.46 & -0.42 & -0.26 \\
-0.46 & 0.33 & -0.26 & -0.11 \\
-0.42 & -0.26 & 2.11 & -0.46 \\
-0.26 & -0.11 & -0.46 & 0.33 \\
\end{array}} \right)
\leq P_S \leq
\left( {\begin{array}{*{20}c}
2.53 & 0.20 & 0.42 & 0.26 \\
0.20 & 0.44 & 0.26 & 0.11 \\
0.42 & 0.26 & 2.53 & 0.20 \\
0.26 & 0.11 & 0.20 & 0.44 \\
\end{array}} \right).
\]
This follows from \textbf{Theorem \ref{thm:period_Q}}. For the Q-piece $\mathcal{Q}_1$ we obtain the following Fenchel Nielsen coordinates $FN_j$ from \textbf{Lemma \ref{thm:FNj}} and the corresponding estimates for $f^u(FN_j)$ and $f^l(FN_j)$:
\begin{table}[htbp]
\begin{center}
\begin{tabular}{|c|c|c|c|c|c|c|}
\hline
$j$ & $\ell(\beta_j)$ & $\ell(\alpha_j)$ & $|{\bf tw}_j|$ & $f^u_{simp}(FN_j)$ & $f^u(FN_j)$ & $f^l(FN_j)$  \\ \hline
1 & 2 & 1 & 0.1 & 0.44 & 0.47 & 0.33 \\ \hline
2 & 2 & 3.032 & 0.017 & 3.16 & 2.53 & 2.11 \\ \hline
12 & 2 & 3.243 & 0.132 & 3.73 & 2.85 & 2.05 \\ \hline
\end{tabular}
\end{center}
\caption{A Q-piece $\mathcal{Q}_1$ given in different Fenchel-Nielsen coordinates $FN_j=(\ell(\beta_j),\ell(\alpha_j),{\bf tw}_j)$ and the values of the corresponding functions $f^u_{simp}(FN_j)$, $f^u(FN_j)$ and $f^l(FN_j)$}
\label{tab:exa}
\end{table}

\section*{Acknowledgement}
The presented work was supported by the Alexander von Humboldt foundation. I would like to thank Peter Buser and Hugo Akrout for helpful discussions and Paman Gujral for proofreading the manuscript. I would also like to thank the referees of the paper for their helpful comments and especially the second referee for his very diligent evaluation of the article.

\end{document}